\documentclass[centertags,leqno,11pt]{article}

\usepackage{amsmath}
\usepackage{amssymb}
\usepackage{latexsym}
\usepackage{amsxtra}
\usepackage{amscd}
\usepackage{theorem}

\usepackage{graphics}
\usepackage{epic}

\theoremstyle{change}

{\theorembodyfont{\slshape}

\newtheorem{thm}{Theorem.}[section]
\newtheorem{cor}[thm]{Corollary.}
\newtheorem{lem}[thm]{Lemma.}
\newtheorem{prop}[thm]{Proposition.}

\newtheorem{defn}[thm]{Definition.}}

{\theorembodyfont{\rmfamily}

\newtheorem{rem}[thm]{Remark.}

}

\renewcommand{\em}{\sl}

\newcommand{\proof}{\noindent {\bf Proof:\ }}
\newcommand{\Endproof}{\hspace*{\fill} $\Box$ \vspace{1ex} \noindent }

\makeatletter
\renewcommand{\subsection}{\@startsection{subsection}{2}%
{\z@}{-3.25ex plus -1ex minus-.2ex}{-1em}{\bf}} \makeatother

\newcommand{\PP}{\mathbb{P}}
\newcommand{\ZZ}{\mathbb{Z}}
\newcommand{\CC}{\mathbb{C}}

\newcommand{\QQ}{\mathbb{Q}}
\newcommand{\NN}{\mathbb{N}}
\newcommand{\FF}{\mathbb{F}}
\renewcommand{\AA}{\mathbb{A}}
\newcommand{\GG}{\mathbb{G}}
\newcommand{\F}{\mathcal{F}}

\newcommand{\V}{\mathcal{V}}

\newcommand{\K}{\mathcal{K}}
\newcommand{\LL}{\mathcal{L}}

\newcommand{\eps}{{\epsilon}}

\newcommand{\bQl}{{\bar{\QQ}_\ell}}
\newcommand{\T}{{\cal T}}
\newcommand{\tame}{{\rm tame}}
\newcommand{\LisseNnl}{{{\rm Lisse}(N,n,\ell)}}
\newcommand{\MT}{{\rm MT}}
\newcommand{\Hyp}{{\rm Hyp}}
\newcommand{\bs}{{\bar{s}}}
\newcommand{\fg}{{\frak{g}}}

\newcommand{\Exc}{{\rm{Exc}}}

\newcommand{\Unip}{{\rm {\bf U}}}

\newcommand{\GL}{{\rm GL}}
\newcommand{\SL}{{\rm SL}}

\newcommand{\Aut}{{\rm Aut}}
\newcommand{\Spec}{{\rm Spec\,}}

\renewcommand{\O}{{\rm O}}
\newcommand{\rig}{{\rm rig}}
\newcommand{\Card}{{\rm Card}}
\newcommand{\HHH}{{\cal H}}
\newcommand{\Ad}{{\rm Ad}}

\newcommand{\MC}{{\rm MC}}

\newcommand{\1}{{\mathbf{1}}}
\newcommand{\ol}{\overline}

\newcommand{\et}{{\rm\acute{e}t}}

\newcommand{\mot}{{\rm mot}}

\newcommand{\To}{\;\longrightarrow\;}

\newcommand{\diag}{{\rm diag}}

\newcommand{\im}{{\rm im}}

\newcommand{\pr}{{{\rm pr}}}

\newcommand{\rk}{{\rm rk}}

\newcommand{\J}{{\rm J}}

\newcommand{\chara}{{\rm char}}

\newcommand{\Corr}{{\rm Corr}}

\newcommand{\Jordan}{{\bf J}}
\newcommand{\FFF}{{\cal F}}\newcommand{\GGG}{{\cal G}}

\newcommand{\G}{{\rm  G}}

\newcommand{\End}{{\rm  End}}

\renewcommand{\H}{{\cal H}}
\renewcommand{\L}{{\cal L}}

\numberwithin{equation}{subsection}
\numberwithin{thm}{subsection}
\theoremstyle{plain}

\begin{document}

\title{Rigid local systems and
 motives of type $G_2$}

\author{Michael Dettweiler and Stefan Reiter}

%\begin{document}
\maketitle 
\begin{abstract} Using the middle convolution functor $\MC_\chi$
which was  
introduced by N.~Katz we 
 prove the existence of 
rigid local systems  whose monodromy
is dense in the simple algebraic group 
$G_2.$ We derive the existence of motives 
for motivated cycles which have a motivic Galois group of type
$G_2.$ 
Granting Grothendieck's standard conjectures,
the existence of motives with motivic Galois group of type $G_2$ can be
deduced, giving a partial answer to a question of Serre. 
\end{abstract}
%\tableofcontents 

%------------------------------------------------------

%\addcontentsline{toc}{section}{Introduction}
%\setcounter{page}{1}
\section*{Introduction}\label{Introduction}

The method of rigidity was first used by B. Riemann 
\cite{Riemann57} in his study of 
 Gau{\ss}' hypergeometric differential equations
${}_2F_1={}_2F_1(a,b,c):$ 
Consider the monodromy representation  
$$\rho: \pi_1^{\rm top}(\PP^1\setminus\{0,1,\infty\},s)\to \GL(V_s)$$
which arises from analytic continuation of the vector space $V_s\simeq \CC^2$ 
of
local solutions of ${}_2F_1$ at $s,$ along paths in $\PP^1\setminus\{0,1,\infty\}$ which are based at $s.$ Let
$\gamma_i,\,i=0,1,\infty,$ 
be simple loops 
around the points $0,1,\infty$ (resp.) which are based at $s.$ 
Then  the monodromy representation 
$\rho$ is {\it rigid} in the sense that 
it is
determined up to isomorphism 
by the Jordan canonical forms of $\rho(\gamma_i),\,i=0,1,\infty.$

One can translate the notion of rigidity into the language of local systems
 by saying that the local system $\L({}_2F_1)$ on 
$\PP^1\setminus\{0,1,\infty\}$ which is given by the
 holomorphic solutions of 
${}_2F_1$  is {\it rigid}
in the following sense:   The {monodromy representation} of $\L({}_2F_1)$
(as defined in \cite{Deligne70}) 
is  determined up to isomorphism by the 
local monodromy representations at the missing points. This 
definition of rigidity extends in the obvious way to other local systems.
Since Riemann's work, the concept of a rigid local system 
has proven to be a very fruitful 
and has 
appeared in many 
different branches of mathematics and physics 
(see e.g. \cite{Beukers-Heckman}, \cite{Ince56}).

A key observation
turned out to be the following: The local sections of the rank-two
 local system
$\L({}_2F_1)$ can be written as linear combinations 
of convolutions $f*g,$ where $f$ and $g$ are solutions of two related
Fuchsian systems of {\it rank one} (see  
to \cite{Katz96}, Introduction).
By interpreting  the convolution as higher direct image and 
using
a transition to 
\'etale sheaves, N. Katz \cite{Katz96}
proved a vast generalization of the above observation:
{ Let  $\FFF$ be  {any} 
irreducible \'etale rigid  local system
on the punctured affine line in the sense specified below. Then 
$\FFF$ can be transformed to a rank one 
sheaf by  a suitable 
iterative application of {\it middle convolutions $\MC_\chi$} and  
tensor products with rank one objects to it (loc.~cit., Chap.~5).} 
(The 
definition of the middle convolution $\MC_\chi$ and its main properties 
 are recalled in Section~\ref{secrevv}.)
This
yields {\it Katz Existence Algorithm} {for irreducible rigid local systems},
which tests, whether a given set of local representations comes from 
an irreducible and rigid local system (loc.~cit., Chap. 6).
This algorithm works simultaneously  in the 
tame \'etale case and in the classical 
case of rigid local systems mentioned above (loc.~cit., Section~6.2 and 6.3).\\

Let $\FFF$ be a 
lisse constructible $\bQl$-sheaf
on a non-empty Zariski open subset $j:U\to \PP^1_{{k}}$ 
which is tamely ramified at the missing 
points $\PP^1_{{k}}\setminus U$ (compare to \cite{SGA5}). 
Call $\FFF$
{\it 
rigid}, if the monodromy representation 
$$\rho_\FFF:\pi_1^\tame(U,\bar{\eta})\To \GL(\FFF_{\bar{\eta}})$$ 
of $\FFF$ is determined 
up to isomorphism 
by the conjugacy classes of the induced representations
of tame inertia groups $I_s^\tame,$ where 
$s\in D:=\PP^1\setminus U.$ Sometimes, we call such a 
sheaf an {\it \'etale rigid  local system}. 
If $\FFF$ is irreducible,
then $\FFF$ is rigid if and only if  the following formula holds:
\begin{equation*}
 \chi(\PP^1,j_*\underline{\End}(\FFF))
=(2-\Card(D))\rk(\FFF)^2+\sum_{s \in D} \dim\left({\rm Centralizer}_{\GL(\FFF_{\bar{\eta}})}(I_s^\tame)\right)=2,
\end{equation*} see \cite{Katz96}, Chap. 2 and 6. 

In preparation to Thm.~1 below, recall that there exist only finitely many 
exceptional simple linear algebraic groups  over an algebraically closed field
which are not isomorphic to a classical group, 
see \cite{Borel}. The smallest of them is the group $G_2$ 
which admits an embedding into the group $\GL_7.$ 
Let us also fix some notation: 
 Let $\1$ (resp. $-\1,$ resp. $\Unip(n)$) denote the 
trivial $\bQl$-valued representation (resp. the unique 
quadratic $\bQl$-valued character, resp. the standard indecomposable 
unipotent $\bQl$-valued representation of degree $n$) of the tame fundamental
 group $\pi_1^\tame(\GG_{m,k}),$ where $k$ is an algebraically
closed field of characteristic $\not=2,\ell.$  The group 
$\pi_1^\tame(\GG_{m,k})$ is isomorphic to the tame inertia 
group $I_s^\tame.$ This can be used to  view representations of 
$I_s^\tame$ as representations of $\pi_1^\tame(\GG_{m,k}).$ 
We prove the following result:\\

\noindent{\bf Theorem~1:} {\it   Let $\ell$ be a prime number and 
let  $k$ be an algebraically closed field of characteristic 
 $\not= 2,\ell.$
Let 
$\varphi,\eta:\pi_1^\tame(\GG_{m,k})\to \bQl^\times $ 
be continuous characters 
such that 
$$ \varphi,\,\eta,\,\varphi\eta,\,\varphi\eta^2,\,\eta\varphi^2,\,
\varphi\ol{\eta} \not\,=\, -\1.$$
Then
there
exists an \'etale rigid  local system
$\HHH(\varphi,\eta)$ of rank $7$ on $\PP^1_k\setminus \{0,1,
\infty\}$ 
whose monodromy group is Zariski dense in
 $G_2(\bQl)$ and whose
local monodromy is as follows: 
\begin{itemize}\item The local monodromy 
at $0$ is of type
$$\,-\1\oplus -\1\oplus -\1\oplus -\1\oplus \1\oplus \1\oplus \1 .$$
\item 
The local monodromy  at $1$ is of type  $$\,\Unip(2)\oplus 
\Unip(2)\oplus \Unip(3).$$
\item The local monodromy  at $\infty$ is of the following form:\end{itemize}
$$ \begin{array}{|c|c|}
\hline 
\mbox{Local monodromy at $\infty$}&\mbox{conditions on $\varphi$ and $\eta$}\\
\hline \hline
\Unip(7)&\varphi=\eta=\1\\
\hline 
\Unip(3,\varphi)
\oplus \Unip(3,\overline{\varphi})\oplus \1& \varphi=\eta\not=\1,\quad 
 \varphi^3=\1\\
\hline 
\Unip(2,\varphi)
\oplus \Unip(2,\ol{\varphi})\oplus 
\Unip(1,\varphi^2)\oplus \Unip(1,\overline{\varphi}^2)
\oplus\1&  \varphi=\eta,\quad  \varphi^4 \not=\1\not= \varphi^6\\
\hline
 \Unip(2,\varphi)
\oplus \Unip(2,\ol{\varphi})\oplus \Unip(3)&  \varphi=\ol{\eta},\quad
\varphi^4\not=  \1\\
\hline
\varphi\oplus \eta \oplus \varphi \eta \oplus 
\overline{\varphi\eta}\oplus \overline{\eta}\oplus \overline{\varphi}\oplus 
\1 &
\varphi, \eta , \varphi\eta , 
\overline{\varphi\eta}, \overline{\eta}, \overline{\varphi}, \1 \\
& \mbox{\rm pairwise different}\\
\hline \end{array}\,\,.$$}

Thm.~1 is proved in a slightly more general form in 
Thm.~\ref{thmneuwicht} below, where it is also proved that
these are the only \'etale rigid local systems of rank $7$ whose monodromy 
is dense in $G_2.$  
The proof of Thm.~\ref{thmneuwicht}
relies heavily 
on Katz' Existence Algorithm.
Using the canonical homomorphism $$\pi_1^{\rm top}(\PP^1(\CC)\setminus 
\{0,1,\infty\})\to \pi_1^\et(\PP^1_\CC\setminus\{0,1,\infty\}),$$
 the existence of a rigid local system
in the classical sense
(corresponding to a representation of the topological 
fundamental group  $\pi_1^{\rm top}(\PP^1(\CC)\setminus 
\{0,1,\infty\})\,$) whose 
monodromy group is Zariski dense in $G_2$ can easily be derived.

Suppose that one has given several 
local representations $$I_s^\tame\To \GL(V)\quad {\rm with}\quad 
 s\in D\cup \{\infty\}$$ 
which are assumed to 
come from an irreducible rigid local system
 $\FFF$ on $\PP^1\setminus D\cup \{\infty\}.$  
It is an empirical observation that the   
 rigidity condition $\chi(\PP^1,j_*\underline{\End}(\FFF))
=2$ and the  (necessary) {\it irreducibility 
condition} $$\chi(\PP^1,j_*\FFF)=(1-\Card(D))\rk(\FFF)
+\sum_{s\in D\cup \{\infty\}} \dim(\FFF_s^{I_s^{\rm tame}})\,\leq \,0$$ 
contradict each other in many cases.
This is  especially often the case if the Zariski closure of the 
monodromy group of $\FFF$ 
is supposed to be small in the underlying 
general linear group.
It is thus astonishing that the above mentioned 
irreducible and rigid 
$G_2$-sheaves exists at all. In fact, the local systems
given by Thm.~1 and Thm.~\ref{thmneuwicht}
 are the first -- and maybe the only -- examples 
of tamely ramified rigid sheaves such that the Zariski closure 
of the monodromy group is an exceptional simple algebraic group. 

We remark
 that in positive 
characteristic, wildly ramified lisse sheaves
on $\GG_m$ with $G_2$-monodromy 
were previously found
by N. Katz (see \cite{KatzGKM}, \cite{Katz90}).
Also, the conjugacy classes in $G_2(\FF_\ell)$ which 
correspond to the local monodromy of the above rigid local system 
$\HHH(\1,\1)$ already occur in the work of Feit, Fong and Thompson 
on the inverse Galois problem (see \cite{FeitFong},
\cite{ThompsonG2}); but only the situation in $G_2(\FF_\ell)$ 
was considered and the transition to rigid local systems was not made.\\

We then apply the above results to give a partial answer to a question
of Serre on the existence of motives with exceptional 
motivic Galois groups. 
 Recall that a  {\it motive} in the Grothendieck sense is 
a triple $M=(X,P,n),\, n\in \ZZ,$ where $X$ is 
a smooth  projective variety over a field $K$ 
 and $P$
is an idempotent correspondence, 
 see e.g. \cite{Saavedra}.
Motives appear in many branches of mathematics (see \cite{Jannsen})
and play a central role in the Langlands program \cite{Langlands}.
Granting Grothendieck's standard conjectures, the category
of Grothendieck motives has the structure of a Tannakian category.
Thus, by the Tannakian formalism,  
 every Grothendieck motive $M$ has conjecturally 
 an algebraic group attached to it,
called the {\it motivic Galois group} of $M$ 
(see \cite{Saavedra} and \cite{DeligneTannaka}).

An unconditional theory of {\em motives for motivated cycles} 
was developed by Andr\'e \cite{Andre96}, who formally adjoins 
a certain homological cycle (the Lefschetz involution) to the algebraic
cycles in order to obtain the Tannakian category of
{ motives for motivated cycles}.  
Let us also mention the 
Tannakian category of 
{\it motives for absolute Hodge cycles}, introduced by Deligne \cite{DMOS}
(for a definition of an absolute Hodge cycle, take a homological 
cycle which satisfies the most visible properties of an algebraic cycle).
In both categories, one has the notion of a motivic Galois
group, given by the Tannakian formalism. It can be shown that any motivated cycle is an absolute Hodge cycle,
so every  motive for motivated cycles is also a motive for absolute Hodge
cycles, see \cite{Andre96}.  Since the category  of motives for motivated cycles
is the minimal extension of Grothendieck's  
category which is unconditionally Tannakian, 
we will work and state our results mainly  in this category. 

The motivic Galois group 
is expected to encode essential properties of a motive. 
Many open conjectures on
motivic Galois groups and related Galois representations are considered
in the article of J.-P. Serre  \cite{Ser2}. 
Under the general 
assumption of Grothendieck's standard conjectures, 
Serre (loc.~cit, 8.8) asks the following question 
{``plus hasardeuse''}:  {\it 
Do there exist motives whose motivic 
Galois group is an exceptional simple algebraic group of type 
$G_2$ (or $E_8$)?} It follows from Deligne's work on 
Shimura varieties that such motives 
cannot be submotives of abelian varieties or 
the motives parametrized by 
 Shimura varieties, see 
\cite{DeligneShimura}. Thus,  
motives with motivic Galois group of type 
$G_2$ or $E_8$ are presumably hard to construct.

There is the notion of a family of motives for motivated cycles,
cf.~\cite{Andre96} and Section~\ref{secresultsfam}.
Using this, we prove the following result (see \ref{thm2new}):\\

\noindent{\bf Theorem~2:} {\it There is
 a  family of motives $M_s$  parametrized  
by  $S= \PP^1\setminus \{0,1,\infty\},$ such that for any  $s \in S(\QQ)$
outside a thin set, the  motive 
$M_s$  has
a motivic Galois group  
of type $G_2.$ }\\

Since the complement of a thin subset of $\QQ$ is infinite (see \cite{SerreMordellWeil}), 
Theorem~2 implies the existence of infinitely many motives for motivated cycles whose motivic Galois 
group is of type $G_2.$
A proof of Thm.~2 will be given in Section~\ref{secmotives}.
It can be shown that under the assumption
of the standard conjectures, the motives $M_s$ are 
Grothendieck motives with motivic Galois group of type $G_2$
(see Rem.~\ref{remext1}).
In this sense, we obtain a positive answer to Serre's question 
in the $G_2$-case. 

The method of the construction of the motives $M_s$ is the motivic 
interpretation of rigid local systems with quasi-unipotent 
local monodromy, introduced by N. Katz
in \cite{Katz96}, Chap. 8.   It follows from Katz' work that 
the sheaf $\HHH(\1,\1)$ in Thm.~1 arises from 
 the cohomology of a smooth { affine} morphism 
$\pi:\Hyp \to \PP^1_{{\QQ}}\setminus \{0,1,\infty\}$
which occurs during   the convolution process 
(see Thm.~\ref{thmmotg2} and Cor. 
\ref{cormotg4}).  
Then a desingularization of 
the relative projective closure of  $\Hyp$
and the work of Y. Andr\'e 
\cite{Andre96} on families of motives imply that a suitable 
compactification and specialization of $\pi$ 
gives  motives over $\QQ$ 
whose motivic Galois groups are of type~$G_2.$ 

In  the Appendix to this paper, which is a joint work 
with N.~ Katz,  the Galois 
representations associated with the above motives $M_s$ are studied. 
It follows from  Thm.~1 of the Appendix., that for two coprime 
integers  $a$ and $b$ which each have at least one 
odd prime divisor,  
the motive $M_{s},\,s=1+\frac{a}{b},$ 
gives rise to $\ell$-adic Galois representations whose 
image is Zariski dense in the group $G_2.$ 
This implies that the motivic Galois group of $M_s$ is of type
$G_2.$  By letting $a$ and $b$ vary among the 
squarefree coprime odd integers $>2,$  one 
obtains  infinitely many non-isomorphic 
motives $M_{1+\frac{a}{b}}$ with motivic Galois group of type $G_2$
(Appendix, Cor.~2~(ii)).

We remark that Gross and Savin \cite{GrossSavin}
propose a completely different 
way to construct motives
with motivic Galois group $G_2$ by looking at 
the cohomology of Shimura varieties of type $G_2$ with { nontrivial 
coefficients.} 
 The connection between these approaches has yet to be 
explored. Due to an observation  of Serre, 
at least the underlying Hodge types coincide.

The authors are indebted to Professors P. Deligne and N. Katz for 
valuable remarks on a previous version of this 
work, and for 
suggesting the generalization of the special 
case of Thm.~1, where $\varphi=\eta=\1,$ to other regular 
characters of the inertia at infinity. We thank 
Prof. J.-P. Serre for his interest in our article and for 
valuable remarks on  the underlying Hodge structures, as well as  
Prof.~Y.~Andr\'e, 
J. Stix,  and  S. Wewers
for helpful comments and discussions. 
Parts of this work where written during the stay 
of the first author at the Institute for Advanced Study in Princeton
in the spring of 2007. He would like to express his thanks for the 
inspirational and friendly atmosphere at the Institute.

\section{Middle convolution and $G_2$-local systems}\label{secrevv}

Throughout the section we fix an algebraically closed field 
$k$ and a prime number $\ell\not= \chara(k).$ 

\subsection{The middle convolution}\label{Secneu1}

Let $G$ be an algebraic group over $k$ and let
$\pi:G\times G\to G$ the multiplication map.
Let $D^b_c(G,\bQl)$ denote the  bounded derived 
category of constructible $\bQl$-sheaves on $G$ (compare to 
\cite{DeligneWeil2}, Section~1, and \cite{Katz96}, Section~2.2).
Given two objects $K,L \in D^b_c(G,\bQl),$ define their {\it $!$-convolution}
 as
$$K*_!L:=R\pi_!(K\boxtimes L)\in D^b_c(G,\bQl)$$ 
and their {\it $*$-convolution} as 
$$K*_*L:=R\pi_*(K\boxtimes L)\in D^b_c(G,\bQl).$$
And element
$K\in D^b_c(G,\bQl)$ 
is called a {\it perverse sheaf} (compare to \cite{BBD}), if 
 $K$  and its dual 
$D(K)$ satisfy $$\dim\left({\rm Supp}(H^i(K))\right)\leq -i,\quad  {\rm resp.} \quad  
\dim\left({\rm Supp}(H^i(D(K)))\right)\leq -i .$$

Suppose that $K$ is a perverse sheaf with the property that for any other
perverse sheaf $L$ on $G,$   the sheaves
$K*_!L$ and $K*_*L$ are again perverse. Then one can define the 
{\it middle convolution} $K*_{\rm mid} L$ of $K$ and $L$  
as the image of $L*_! K$ in $ L*_* K$ under the ``forget supports map''
in the abelian category of perverse sheaves.

Let us now consider the situation, where 
$G=\AA^1_k:$ 
For any nontrivial continuous character
$$ \chi: \pi_1^\tame(\GG_{m,k})\To \bQl^\times,$$ let $\LL_\chi$ denote the 
corresponding lisse
 sheaf of rank one on $\GG_{m,k}.$ Let $j:\GG_m\to \AA^1$
denote the inclusion. 
From $j_*\L_\chi$ one obtains
a perverse sheaf $j_*\L_\chi[1]$ 
on $\AA^1$ by placing the sheaf in degree $-1.$ 
Since $!$-convolution (resp. $*$-convolution) with 
$j_*\L_\chi[1]$ preserves perversity (see \cite{Katz96}, Chap. 2), the 
middle 
convolution $K*_{\rm mid} j_*\L_\chi[1]$ 
is defined for any perverse sheaf $K$ on $\AA^1.$ \\

The following notation will be used below: For any scheme
$W$ and any map $f: W\to \GG_m,$ define 
\begin{equation}\label{eqkatzrem} \L_{\chi(f)}:=f^*\L_{\chi}.\end{equation} 
The identity character will be denoted by 
$\1$ and $-\1$ denotes the unique quadratic 
character of $\pi_1^\tame(\GG_m).$ The {\it inverse 
character} of $\chi$ will be denoted by $\overline{\chi}$
(by definition, $\chi\otimes \ol{\chi}=\1$).\\

The following category will be of importance below: 

\begin{defn} {\rm Let 
 $\T_\ell=\T_\ell(k)$ denote the full subcategory of constructible 
$\bQl$-sheaves $\FFF$ on $\AA^1_k$ 
which satisfy the following conditions: 
\begin{itemize}
\item There exists a dense open subset $j:U\to \AA^1$ such that 
$j^*\FFF$ is lisse and irreducible on $U,$ and such that 
$\FFF\simeq j_*j^*\FFF.$
\item The lisse sheaf $j^*\FFF$ is tamely ramified at every point 
of $\PP^1\setminus U.$
\item There are at least two distinct points of $\AA^1$ at which $\FFF$
fails to be lisse.
\end{itemize}}
\end{defn}

 The properties of  $\T_\ell$ imply that
 $\FFF[1]*_{\rm mid} \L_\chi[1]$ is a single sheaf placed in degree $-1$
(\cite{Katz96}, Chap. 5), leading to the 
 {\it middle convolution functor}
$$ \MC_\chi: \T_\ell \to \T_\ell,\quad \FFF \mapsto (\FFF[1]*_{\rm mid} \L_\chi[1])[-1],$$ see \cite{Katz96}, 5.1.5. (Note that, by the 
definition of $\T_\ell,$ the sheaf $\MC_\chi(\FFF)$ 
is again irreducible, compare to loc.~cit. Thm.~3.3.3) \\

An important property of $\MC_\chi$ is the following:
\begin{equation}\label{eqnneu9} \MC_\chi\circ \MC_\rho=\MC_{\chi\rho}\quad \mbox{\rm if}\quad 
\chi\rho\neq \1,\quad 
{\rm and} \quad   \MC_\chi\circ \MC_{\ol{\chi}}={\rm Id}.\end{equation}
Let $U\subseteq \AA^1$ be an open subset of $\AA^1$ such that 
$\FFF|_U$ is lisse and let $\iota:U\to \PP^1$ be the canonical inclusion.
The sheaf $\FFF\in \T_\ell$ is called {\it cohomologically rigid}, if
the {\it index of rigidity} 
$$ \rig(\FFF)=\chi(\PP^1,\iota_*(\underline{\End}(\FFF|_{U})))
$$ is equal to $2.$ Then $\MC_\chi$ carries rigid 
elements in $\T_\ell$ to rigid elements in $\T_\ell$ by the following result,
see \cite{Katz96}, 6.0.17:
\begin{equation}\label{eqnneu8}
\rig(\FFF)=\rig(\MC_\chi(\FFF)).\end{equation}

\subsection{The numerology of the middle convolution}
We recall the effect of the middle convolution 
on the Jordan canonical forms of the 
local monodromy, given by Katz in \cite{Katz96}, Chap. 6:\\

Let $\FFF\in \T_\ell$ and let $j:U\mapsto \AA^1_x$ denote an open subset
such that $j^*\FFF$ is lisse. Let $D:=\AA^1\setminus U.$
 Then, for any point $s\in D \cup \{\infty\}=\PP^1\setminus U,$ the sheaf 
$\FFF$ gives rise to the {\it local monodromy representation $\FFF(s)$}
of the 
tame inertia subgroup $I(s)^\tame$  (of the absolute 
Galois group of the generic point of $\AA^1$) at $s.$ The representation
$\FFF(s)$ decomposes as follows as a direct sum 
(character)$\otimes$(unipotent representation), 
where the sum is over the set of 
continuous $\bQl$-characters $\rho$ of 
$\pi_1^\tame(\GG_{m,k})\simeq I(s)^\tame:$
$$ \FFF(s)=\bigoplus_{\rho} \L_{\rho(x-s)}\otimes {\rm Unip}(s,\rho,\FFF)\quad 
\mbox{\rm for all} \quad s\in D$$
and 
$$ \FFF(\infty)=\bigoplus_{\rho} \L_{\rho(x)}
\otimes {\rm Unip}(\infty,\rho,\FFF).$$
Here, the following convention is used: If one starts with a rank 
one object $\FFF,$ which at $s\in D$ 
gives locally rise to a character 
$\chi_s$ of $\pi_1^\tame(\GG_m)$ then 
$$\chi_\infty=\prod_{s\in D} \chi_s.$$

For $s\in D\cup \{\infty\},$ write 
 ${\rm Unip}(s,\rho,\FFF)$ as a direct sum of Jordan blocks of lengths
$\{n_i(s,\rho,\FFF)\}_i.$ This leads to a decreasing sequence of non-negative
integers
$$ e_1(s,\rho,\FFF)\geq e_2( s,\rho,\FFF)\geq \ldots\geq 
e_k( s,\rho,\FFF)=0\quad {\rm for} \quad k \gg 0,$$
where the number $e_j( s,\rho,\FFF)$ is defined to be the number of 
Jordan blocks in ${\rm Unip}(s,\rho,\FFF)$ whose length is $\geq j.$ 

\begin{prop}\label{propnum} Let $\FFF\in \T_\ell$ be of generic rank $n.$
Then the following holds:
\begin{enumerate}
\item 
\begin{eqnarray*}\rk(\MC_\chi(\FFF))&=&\sum_{s\in D}\rk(\FFF(s)/
(\FFF(s)^{I(s)}))-
\rk((\FFF(\infty)\otimes \L_{{\chi}})^{I(\infty)})\\
&=&\sum_{s\in D} (n-e_1(s,\1,\FFF))-e_1(\infty,\ol{\chi},\FFF).\end{eqnarray*}
\item
For $s\in D$ and $i\geq 1,$  the following holds:
\begin{eqnarray*}
 e_i(s, \rho\chi, \MC_\chi(\FFF))&=&e_i(s,\rho,\FFF)\quad {\rm if} \quad
\rho\not= \1 \quad {\rm and}\quad  \rho\chi\not= \1,\\
 e_{i+1}(s,\1,\MC_\chi(\FFF))&=&e_i(s,\bar{\chi},\FFF),\\
 e_i(s,\chi,\MC_\chi(\FFF))&=&e_{i+1}(s,\1,\FFF). \end{eqnarray*}
Moreover, 
$$e_1(s,\1,\MC_\chi(\FFF))=\rk(\MC_\chi(\FFF))-n+e_1(s,\1,\FFF).$$

\item For $s=\infty$ and $i\geq 1,$ the following holds:
\begin{eqnarray*}
e_i(\infty, \rho{\chi}, \MC_\chi(\FFF))&=&e_i(\infty,\rho,\FFF)\quad {\rm if} \quad
\rho\not= \1\quad {\rm and}\quad  \rho\chi\not= \1,\\
e_{i+1}(\infty,\chi,\MC_\chi(\FFF))&=&e_i(\infty,\1,\FFF),\\
e_i(\infty,\1,\MC_\chi(\FFF))&=&e_{i+1}(\infty,\ol{\chi},\FFF).\end{eqnarray*}
Moreover, 
$$e_1(\infty,{\chi},\MC_\chi(\FFF))=
\sum_{s\in D}(\rk (\FFF)-e_1(s,\1,\FFF)) -\rk(\FFF).$$
\end{enumerate}
\end{prop}

\proof Claim (i) is \cite{Katz96}, Cor. 3.3.7. 
The first three equalities in (ii) are loc.~cit., 6.0.13. 
The last equality in (ii) follows from loc.~cit., 6.0.14.
To deduce (iii), we argue as follows: From loc.~cit., 3.3.6, and 6.0.5,
for any $\FFF\in \T_\ell$ there exists an $I(\infty)^\tame$-representation $M(\infty,\FFF)$ of 
rank $\sum_{s\in D}(n-e_1(s,\1,\FFF))$ with the 
following  properties: 
\begin{eqnarray*} 
E_i(\infty,\rho,\FFF)&=
&e_i(\infty,\rho,\FFF)\quad {\rm if}\quad \rho\not= \1,\\
E_{i+1}(\infty,\1,\FFF)&=&e_i(\infty,\1,\FFF)\quad {\rm for} \quad i\geq 1,\\
E_1(\infty,\1,\FFF)&=&\rk(M(\infty,\FFF))-\rk(\FFF),\end{eqnarray*} where 
the numbers $E_i(\infty,\rho,\FFF)$ denote the invariants
associated to $M(\infty,\FFF)$ which are defined analogously as the invariants 
$e_i(s,\rho,\FFF)$ for $\FFF(s).$
Moreover, by loc.~cit., 6.0.11, the following holds:
$$ E_i(\infty,\rho\chi,\MC_\chi(\FFF))=E_i(\infty,\rho,\FFF),\quad \mbox{\rm for all}\quad
i\geq 1 \quad \mbox{\rm and } \rho.$$
By combining the last equations, it follows that if $\rho\chi\not= \1$ and $\rho\not= \1,$ then 
$$ e_i(\infty,\rho\chi,\MC_\chi(\FFF))=E_i(\infty,\rho\chi,\MC_\chi(\FFF))=
E_i(\infty,\rho,\FFF)=e_i(\infty,\rho,\FFF).$$
If  $\rho=\1,$ since $\chi$ is nontrivial, the following holds 
$$ e_{i+1}(\infty,\chi,\MC_\chi(\FFF))=E_{i+1}(\infty,\chi,\MC_\chi(\FFF))=
E_{i+1}(\infty,\1,\FFF)=e_i(\infty,\rho,\FFF).$$
Moreover,
$$ e_i(\infty,\1,\MC_\chi(\FFF))=E_{i+1}(\infty,\1,\MC_\chi(\FFF))=E_{i+1}(\infty,\ol{\chi},\FFF)
=e_{i+1}(\infty,\ol{\chi},\FFF),$$ since $\chi$ and thus $\ol{\chi}$ are nontrivial. Finally, 
\begin{eqnarray*}
 e_1(\infty,\chi,\MC_\chi(\FFF))
&=&E_1(\infty, \chi,\MC_\chi(\FFF))\\
&=&E_1(\infty,\1,\FFF)\\
&=&
\rk(M(\infty,\FFF))-\rk(\FFF)\\
&=&\sum_{s\in D}(\rk (\FFF)-e_1(s,\1,\FFF)) -\rk(\FFF),\end{eqnarray*}
where the last equality follows from 
 loc.~cit., 6.0.6.
\Endproof

Let $\FFF\in \T_\ell$ and let $\L$ be a middle extension sheaf on $\AA^1$ (i.e., there 
exists an open subset $j:U\to \AA^1$ such that $j^*\L$ is lisse  and such that 
$\L\simeq j_*j^*\L).$ Assume that $\FFF|_U$ is also lisse. Then 
the {\it middle tensor product} of $\FFF$ and $\L$ is defined as 
$$\MT_\L(\FFF)=j_*(\FFF|_U\otimes \L|_U),$$
compare to loc.~cit., 5.1.9. Obviously, the generic rank of 
$\MT_\L(\FFF)$ is the same as the generic rank of $\FFF.$ 
For any $s\in D\cup\{\infty\},$ denote by 
$\chi_{s,\L}$ the unique character $\rho$ with $e_1(s,\rho,\L)=1.$
Then the following holds (loc.~cit., 6.0.10): 
\begin{equation}\label{eqneu1}
e_i(s,\rho\chi_{s,\L},\MT_\L(\FFF))=e_i(s,\rho,\FFF).\end{equation}

\subsection{Classification of irreducible rigid local systems with 
$G_2$-monodromy}\label{secwichtig}

In this section, we give a complete classification of 
rank $7$ rigid sheaves $\HHH\in \T_\ell$ 
whose associated monodromy
group is Zariski dense in the group $G_2(\bQl).$  

Let us first  collect the information on the
conjugacy classes of the simple algebraic group $G_2$ which are needed
below. 
In Table~1 below, we
lists the possible Jordan canonical forms 
of elements of the group $G_2(\bQl)\leq \GL_7(\bQl)$ together with the 
dimensions of the centralizers in the group $G_2$ and in the 
group $\GL_7(\bQl).$ That the list exhausts all possible cases can be seen 
 (using Jordan decomposition) from the fact
that a semisimple element in $G_2(\bQl)\leq \GL_7(\bQl)$ is of the form 
$ \diag(x,y,xy,1,(xy)^{-1},y^{-1},x^{-1}),$  and from  
the same arguments as in  
\cite{Lawther} for the classification 
of the unipotent classes. The computation of the centralizer dimension 
in $G_2(\bQl)$ follows 
using the same arguments as in \cite{Chang}, the centralizer dimension 
in the group $\GL_7$ follows e.g. from \cite{Katz96}, 3.1.15.

We use the following conventions: $E_n\in \bQl^{n\times n}$
 denotes the identity matrix, 
$\Jordan(n)$ denotes a unipotent 
Jordan 
block of length $n,$ $\epsilon\in \bQl^\times$ 
 denotes a primitive $3$-rd root
of unity,  and $i\in \bQl^\times$ denotes a primitive $4$-th root of unity.
Moreover, an expression like $(x\Jordan(2),x^{-1}\Jordan(2),x^2,x^{-2},1)$
denotes a matrix in Jordan canonical form 
in $\GL_7(\bQl)$ with one Jordan block of 
length 2 having eigenvalue $x,$ one Jordan block of 
length 2 having eigenvalue $x^{-1},$ and three Jordan blocks of length $1$ 
having eigenvalues $x^2,x^{-2},1$ (resp.). \\

\begin{equation*}\label{table1}
\begin{array}{|c|c|c|}
\hline
 \textrm{Jordan form}  &\mbox{\rm Centralizer dimension in} &  
\textrm{Conditions}\\
               & G_2  \quad\quad\quad \GL_7&  \\
\hline 
 E_7& 14 \quad\quad\quad  49&\\
(\Jordan(2),\Jordan(2),E_3)& 8\quad\quad\quad 29&\\
 (\Jordan(3),\Jordan(2),\Jordan(2))&6\quad\quad\quad 19&\\
(\Jordan(3),\Jordan(3),\1)& 4\quad\quad\quad 17&\\
 \Jordan(7)&2 \quad\quad\quad 7&\\
\hline
 (-E_4,E_3)& 6\quad\quad\quad  25 & \\
(-\Jordan(2),-\Jordan(2),E_3)&4\quad\quad\quad 17&\\
(-\Jordan(2),-\Jordan(2),\Jordan(3))&4\quad\quad\quad 11&\\
(-\Jordan(3),-\1,\Jordan(3))&2\quad\quad\quad 9&\\
\hline\end{array}
\end{equation*}

\begin{equation*}
\begin{array}{|c|c|c|}
\hline
 \textrm{Jordan form}  &\mbox{\rm Centralizer dimension in} &  
\textrm{Conditions}\\
               & G_2  \quad\quad\quad \GL_7&  \\
\hline 
(\eps E_3,1,\eps^{-1}E_3 )& 8\quad\quad\quad 19 & \\
(\eps \Jordan(2),\eps^{-1}\Jordan(2),\eps,\eps^{-1},1)&4\quad\quad\quad 11&\\
(\eps \Jordan(3),\eps^{-1}\Jordan(3),1 )&2\quad\quad\quad 7&\\
\hline 
(i,i,-1,1,i^{-1},i^{-1},-1) &4\quad\quad\quad   13&\\
(i\Jordan(2),i^{-1}\Jordan(2),-1,-1,1)&2\quad\quad\quad 9 &\\
\hline
(x,x,x^{-1},x^{-1},1,1,1)&4 \quad\quad\quad 17& x^2\not=1\\
(x,x,x^2,1,x^{-1},x^{-1},x^{-2})&4 \quad\quad\quad 11&x^4\not=1\not=x^3\\
(x,-1,-x,1,-x^{-1},-1,x^{-1})&2\quad\quad\quad 9& x^4\not=1 \\
(x\Jordan(2),x^{-1}\Jordan(2),x^2,x^{-2},1)&2\quad\quad\quad 7 & x^4\not=1\\
(x\Jordan(2),x^{-1}\Jordan(2),\Jordan(3))&2 \quad\quad\quad 7& x^2\not=1\\
\hline
(x,y,xy,1,(xy)^{-1},y^{-1},x^{-1})&2\quad\quad\quad 7 &x,y,xy,1,(xy)^{-1},y^{-1},x^{-1}\\
&& \textrm{pairwise different}\\
\hline
\end{array}
\end{equation*}
\vspace{.05cm}
\begin{center} {\bf Table~1:} The $\GL_7$ conjugacy classes of $G_2$
\end{center}
\vspace{.1cm}

We will also use the following notation in Thm.~\ref{thmneuwicht} below:  
 Let 
$\Unip(i)$ denote the  $\bQl$-valued 
representation of $\pi_1^\tame(\GG_m)$ which sends a 
generator of $\pi_1^\tame(\GG_m)$ to the Jordan block 
$\J(i).$ For any character $\chi$ of 
$\pi_1^\tame(\GG_m)$ let  
$$\Unip(i,\chi):=\chi\otimes \Unip(i),$$ let 
$$ -\Unip(i):=-\1\otimes \Unip(i)$$ 
($-\1$ denoting the unique quadratic character of  $\pi_1^\tame(\GG_m)$),
and let 
$\Unip(i,\chi)^j$ denote the $j$-fold direct 
sum of the representation $\Unip(i,\chi).$ 
%The main result of this section is the following:

\begin{thm} \label{thmneuwicht} 
Let $\ell$ be a prime number and 
let  $k$ be an algebraically closed field with $\chara(k)\not= 2,\ell.$
 Then the following holds:
\begin{enumerate}
\item  
Let 
$\alpha_1,\alpha_2\in \AA^1(k)$ be two 
disjoint points and let
$\varphi,\eta:\pi_1^\tame(\GG_{m,k})\to \bQl^\times $ 
be continuous characters 
such that 
$$ \varphi,\,\eta,\,\varphi\eta,\,\varphi\eta^2,\,\eta\varphi^2,\,
\varphi\ol{\eta} \not\,=\, -\1.$$
Then
there
exists an irreducible  cohomologically rigid  
sheaf $\HHH=\HHH(\varphi,\eta)\in \T_\ell(k)$ of generic rank $7$
whose local monodromy  is as follows: 
\begin{itemize}\item The local monodromy 
at $\alpha_1$ is
$\,-\1^4\oplus \1^3 .$
\item 
The local monodromy  at $\alpha_2$ is  $\,\Unip(2)^2\oplus \Unip(3).$
\item The local monodromy   at $\infty$ is of the following form:
$$ \begin{array}{|c|c|}
\hline 
\mbox{Local monodromy at $\infty$}&\mbox{conditions on $\varphi$ and $\eta$}\\
\hline \hline
\Unip(7)&\varphi=\eta=\1\\
\hline 
\Unip(3,\varphi)
\oplus \Unip(3,\overline{\varphi})\oplus \1& \varphi=\eta\not=\1,\quad 
 \varphi^3=\1\\
\hline 
\Unip(2,\varphi)
\oplus \Unip(2,\ol{\varphi})\oplus 
\Unip(1,\varphi^2)\oplus \Unip(1,\overline{\varphi}^2)
\oplus\1&  \varphi=\eta,\quad  \varphi^4 \not=\1\not= \varphi^6\\
\hline
 \Unip(2,\varphi)
\oplus \Unip(2,\ol{\varphi})\oplus \Unip(3)&  \varphi=\ol{\eta},\quad
\varphi^4\not=  \1\\
\hline
\varphi\oplus \eta \oplus \varphi \eta \oplus 
\overline{\varphi\eta}\oplus \overline{\eta}\oplus \overline{\varphi}\oplus 
\1 &
\varphi, \eta , \varphi\eta , 
\overline{\varphi\eta}, \overline{\eta}, \overline{\varphi}, \1 \\
& \mbox{\rm pairwise different}\\
\hline \end{array}$$
 \end{itemize}
Moreover,  the restriction
$\HHH|_{\AA^1_k\setminus\{\alpha_1,\alpha_2\}}$ is lisse
and the monodromy group 
associated to $\HHH$ is a Zariski dense subgroup of 
the simple exceptional  algebraic group $G_2(\bQl).$ 

\item Assume that $\HHH\in \T_\ell$
is a cohomologically rigid $\bQl$-sheaf of generic rank $7$ 
which fails to be lisse at $\infty$ and 
such that the monodromy group 
associated to $\HHH$ is Zariski dense in the group $G_2(\bQl).$
Then $\HHH$ fails to be lisse
at exactly two disjoint points $\alpha_1,\alpha_2\in \AA^1(k)$ and,
up to a permutation of the points $\alpha_1,\alpha_2,\infty,$ 
the above  list exhaust all the possible 
local monodromies of $\HHH.$ 
\end{enumerate}
\end{thm}

\proof Let us introduce the following notation:
Let
$j:U:=\AA^1_k\setminus \{\alpha_1,\alpha_2\}\to \PP^1$ denote the tautological inclusion. Let $x-\alpha_i,\,i=1,2,$ denote the morphism 
$U\to \GG_m$ which is induced by sending $x\in U$ to $x-\alpha_i.$ 
For any pair of continuous
characters $\chi_1,\chi_2:\pi_1^\tame(\GG_m)\to \bQl^\times,$ 
let $$\L(\chi_1,\chi_2):=j_*\left(\L_{\chi_1(x-\alpha_1)}\otimes 
\L_{\chi_2(x-\alpha_2)}\right)$$ 
(using the notation of Formula~\ref{eqkatzrem}). 
Let 
$$\F_1=\L(-\1,-\varphi\eta)\in \T_\ell.$$
Define inductively a sequence of sheaves 
$\HHH_0,\ldots,\HHH_6$ in $\T_\ell,$
by setting 
$$ \HHH_0:=\F_1,\quad {\rm and}\quad \HHH_{i}:=\MT_{\F_{i+1}}(\MC_{\rho_i}(\HHH_{i-1})),
\quad {\rm for}\quad i=1,\ldots,6,$$
where the $\F_i$ and $\rho_i$ are as follows:
$$ \F_3=\F_5=\F_7=\L(-\1,\1),\quad  \F_2=\L(\1,-\overline{\varphi}),\quad \F_4=\L(\1,-\varphi\overline{\eta}),\quad
\F_6=\L(\1,-\overline{\varphi})$$
and $$ \rho_1:=-\overline{\varphi\eta^{2}},\; \rho_2:=-\varphi\eta^2,\;
\rho_3:=-\overline{\varphi\eta},\; \rho_4:=-\varphi\eta,\; \rho_5:=-\ol{\varphi},\; 
\rho_6:=-\varphi.$$

We now distinguish $5$ cases, which correspond to the 
different types of local monodromy  at $\infty$ 
listed above:

\noindent{\it Case 1:} Let $\varphi=\eta=\1.$
 The following table lists the local monodromies 
of the sheaves $\HHH_0,\ldots,\HHH_6=\HHH$ at the points 
$\alpha_1,\alpha_2,\infty$  (the proof is a direct computation, using
Prop.~\ref{propnum} and Equation \eqref{eqneu1}):
$$ \begin{array}{|c|c|c|c|}
\hline 
&{\rm at}\quad \alpha_1&{\rm at}\quad \alpha_2&{\rm at}\quad \infty\\
\hline
\HHH_0& -\1&-\1&\1\\
\hline 
\HHH_1&  \Unip(2)& -\Unip(2)& \Unip(2)\\
\hline 
\HHH_2&  -\1^2 \oplus \1& \Unip(3)& \Unip(3)\\
\hline
\HHH_3& \Unip(2)^2& \Unip(2)\oplus -\1^2
& \Unip(4)\\
\hline 
\HHH_4& \1^2\oplus -\1^3&
 \Unip(2)^2\oplus -\1&\Unip(5)\\
\hline 
\HHH_5 &\Unip(2)^3 & 
-\Unip(2)\oplus-\1^2\oplus\1^2 & \Unip(6)\\
\hline
\HHH_6 & -\1^4\oplus \1^3&
\Unip(2)^2\oplus\Unip(3)&\Unip(7) \\
\hline
\end{array}$$
By  the results of Section~\ref{Secneu1} and Prop.~\ref{propnum}, the sheaf 
$\HHH=\HHH_6$ is a cohomologically rigid 
 irreducible sheaf of rank $7$ 
in $\T_\ell$ which is 
lisse on the open subset
$U=\AA^1_k\setminus\{\alpha_1,\alpha_2\}\subseteq \AA^1_k.$ 
The lisse sheaf $\HHH|_U$ corresponds to 
a representation 
$$\rho: \pi_1^\tame(\AA^1\setminus \{\alpha_1,\alpha_2\})\to \GL(V),$$
where $V$ is a $\bQl$-vector space of dimension $7.$ Let $G$ be the image of $\rho.$ 
Note that $G$ is an irreducible subgroup of $\GL(V)$ since 
$\HHH$ is irreducible.
In the following, we fix an isomorphism $V\simeq \bQl^7.$ This induces
 an isomorphism $\GL(V)\simeq \GL_7(\bQl),$ so we can view
$G$ as a subgroup of $\GL_7(\bQl).$  

We want to show that $G$ is contained in a conjugate of 
the group $G_2(\bQl)\leq \GL_7(\bQl).$ For this we argue as 
in \cite{Katz90}, Chap. 4.1: First note that the 
local monodromy at $s\in \{\alpha_1,\alpha_2,\infty\}$ 
 can be (locally) conjugated in $\GL_7(\bQl)$ 
into the orthogonal group
$O_7(\bQl).$ 
It follows thus from the rigidity of the representation 
 $\rho$ that  
there exists an element $x\in \GL(V)$ such that 
$$ {\rm Transpose}(\rho(g)^{-1})\,=\,\rho(g)^x\quad \forall g \in 
\pi_1^\tame(\AA^1\setminus \{\alpha_1,\alpha_2\}).$$ 
In other words, the group $G$ respects the nondegenerate
bilinear form 
given by the element $x^{-1}.$ Since $G$ is irreducible and the 
dimension of $V$ is $7,$ this form 
has to be symmetric. Thus we can assume that  $G$ 
is contained in the orthogonal group $O_7(\bQl).$
  By the results of Aschbacher
\cite{Asch}, Thm.~5~(2)~and~(5), an irreducible 
subgroup $G$ of $\O_7(K)$ ($K$ denoting  an 
algebraically closed field or a finite field)
  lies inside an $\O_7(K)$-conjugate 
of $G_2(K),$ if and only if $G$  has a nonzero invariant
in the third exterior power $\Lambda^3(V)$ of $V=K^7.$ 
In our case, this  is equivalent to 
\begin{equation}\label{eqqq}
H^0(U,\Lambda^3(\tilde{\HHH}))\simeq H^0(\PP^1,j_*\Lambda^3(\tilde{\HHH}))\not= \{0\},\end{equation} where $\tilde{\HHH}=\HHH|_U.$ 
Poincar\'e duality
implies that the last formula is equivalent to  
\begin{equation}\label{eqqq1}H^2_c(U,\Lambda^3(\tilde{\HHH}))\simeq H^2(\PP^1,j_*\Lambda^3(\tilde{\HHH}))\not= \{0\}.\end{equation}
 The Euler-Poincar\'e
Formula implies that 
\begin{eqnarray}\label{eqabove}
\quad \quad \,\,\,\,\chi(\PP^1,j_*\Lambda^3(\tilde{\HHH}))& =& 
h^0(\PP^1,j_*\Lambda^3(\tilde{\HHH}))- h^1(\PP^1,j_*\Lambda^3(\tilde{\HHH}))
+ h^2(\PP^1,j_*\Lambda^3(\tilde{\HHH}))\\
& =&
\chi(U)\cdot \rk(\Lambda^3(\tilde{\HHH}))
+\sum_{s\in \{\alpha_1,\alpha_2,\infty\}} 
\dim(\Lambda^3(\tilde{\HHH})^{I(s)})\nonumber\\
&=& -35 +19+13+5\nonumber \\
&=& 2\nonumber\end{eqnarray}
(note that $\chi(U)=h^0(U)-h^1(U)+h^2(U)=-1,$ $\rk(\Lambda^3(\tilde{\HHH}))=35,$ 
and that for 
 $s=\alpha_1,\alpha_2,\infty,$ the dimension of the local invariants  
$\dim(\Lambda^3(\tilde{\HHH})^{I(s)})$ can be computed to 
be  equal to $19,13,5,$ resp.).
 It follows thus from the equivalence of  \eqref{eqqq} and \eqref{eqqq1}
that 
$h^0(U,\Lambda^3(\tilde{\HHH}))\geq 1.$ Therefore,
the monodromy group $G$ can be assumed to 
be contained in $G_2(\bQl).$ Let $\ol{G}$ denote the Zariski closure of 
$G$ in $G_2(\bQl).$ 
By \cite{Asch}, Cor. 12, a Zariski closed proper maximal 
subgroup of $G_2(\bQl)$ is either 
reducible or $G$ is isomorphic to the group  $\SL_2(\bQl)$ acting on the 
vector space of homogeneous polynomials of degree $6$. 
In the latter case, the unipotent elements of the image of 
$\SL_2(\bQl)$ are equal to the identity matrix, or they are conjugate 
in $\GL_7(\bQl)$ to a Jordan 
block of length $7.$ 
Since the local monodromy of $\HHH$ at $\alpha_2$ is not of this form, 
$G$ must coincide with $\G_2(\bQl).$
This finishes the proof of Claim (i) in the Case 1.

In the following, we only list the local monodromy of the of 
the sheaves $$\HHH_0,\ldots,\HHH_6=\HHH$$ in the remaining
Cases 2-5. In each case, rigidity implies that the image of $\pi_1(U)$ 
is contained in an orthogonal group, and one can compute that  an analogue
of Formula \eqref{eqabove} holds. Thus the image of $\pi_1(U)$ 
is Zariski dense in $G_2$ by the same arguments as in Case 1. 

\noindent{\it Case 2:}  $\varphi=\eta$ and $\varphi$ is  
nontrivial of order $3.$ Then 
$$ \begin{array}{|c|c|c|c|}
\hline 
&{\rm at}\quad \alpha_1&{\rm at}\quad \alpha_2 &{\rm at}\quad \infty\\
\hline
\HHH_0& -\1&-\ol{\varphi}&\ol{\varphi}\\
\hline 
\HHH_1&  \Unip(2)& -\varphi\oplus -\ol{\varphi}& \varphi\oplus \ol{\varphi}\\
\hline 
\HHH_2&  -\1^2 \oplus \1& \varphi\oplus \ol{\varphi}\oplus \1&
 \varphi\oplus\ol{\varphi}\oplus \1\\
\hline
\HHH_3& \Unip(1,\varphi)^2\oplus \1^2&
\ol{\varphi}\oplus \1 \oplus \Unip(1,-\1)^2&
 \Unip(2,\varphi)\oplus \1\oplus \ol{\varphi}\\
\hline 
\HHH_4& \1^2\oplus -\1^3&
-\varphi\oplus \Unip(1,\ol{\varphi})^2\oplus \1^2&
\Unip(2)\oplus \Unip(2,\ol{\varphi})\oplus \varphi\\
\hline 
\HHH_5 &\Unip(1,\ol{\varphi})^3\oplus \1^3& 
\Unip(2,-\varphi)\oplus \1^2\oplus \Unip(1,-\varphi)^2&
\Unip(3,\varphi)\oplus \Unip(2)\oplus \ol{\varphi}\\
\hline
\HHH_6 & -\1^4\oplus \1^3&
\Unip(2)^2\oplus\Unip(3)&\Unip(3,\varphi)\oplus \Unip(3,\ol{\varphi})\oplus 
\1 \\
\hline
\end{array}.$$

\noindent{\it Case 3:}  $\varphi=\eta$ and  
$\varphi^4\not=\1\not=\varphi^6.$ Then 

$$ \begin{array}{|c|c|c|c|}
\hline 
&{\rm at}\quad \alpha_1&{\rm at}\quad \alpha_2&{\rm at}\quad \infty\\
\hline
\HHH_0& -\1&-{\varphi}^2&{\varphi}^2\\
\hline 
\HHH_1&  \ol{\varphi}^3\oplus \1&
 -\ol{\varphi}^2\oplus -\ol{\varphi}& \ol{\varphi}^4\oplus \ol{\varphi}^2\\
\hline 
\HHH_2&  -\1^2 \oplus \1&
 \varphi\oplus {\varphi}^2\oplus \1&
 \varphi^3\oplus\ol{\varphi}\oplus \varphi\\
\hline
\HHH_3& \Unip(1,\ol{\varphi})^2\oplus \1^2&
\ol{\varphi}\oplus \1 \oplus \Unip(1,-\1)^2&
\ol{\varphi}^2\oplus \varphi\oplus \ol{\varphi}^3\oplus \ol{\varphi}\\
\hline 
\HHH_4& \1^2\oplus -\1^3&
-\varphi\oplus \Unip(1,{\varphi})^2\oplus \1^2&
\varphi^2\oplus \1\oplus \varphi^3\oplus \ol{\varphi}\oplus \varphi\\
\hline 
\HHH_5 &
\Unip(1,\ol{\varphi})^3\oplus \1^3& 
\Unip(2,-\ol{\varphi})\oplus \1^2\oplus \Unip(1,-\ol{\varphi})^2&
\1\oplus \Unip(2,\ol{\varphi}^2)\oplus \varphi\oplus \ol{\varphi}^3\oplus 
\ol{\varphi}\\
\hline
\HHH_6 & -\1^4\oplus \1^3&
\Unip(2)^2\oplus\Unip(3)&\Unip(2,\varphi)\oplus \Unip(2,\ol{\varphi})
\oplus \varphi^2\oplus \ol{\varphi}^2\oplus \1 \\
\hline
\end{array}$$

\noindent{\it Case 4:}  $\varphi=\ol{\eta}$ and 
$\varphi^4\not=\1.$ Then 
$$ \begin{array}{|c|c|c|c|}
\hline 
&{\rm at}\quad \alpha_1&{\rm at}\quad \alpha_2&{\rm at}\quad \infty\\
\hline
\HHH_0& -\1&-\1&\1 \\
\hline 
\HHH_1&  {\varphi}\oplus \1&
 -\1\oplus -\ol{\varphi} &
\Unip(2)\\
\hline 
\HHH_2&  -\1^2 \oplus \1&
 \ol{\varphi}\oplus \ol{\varphi}^2\oplus \1&
 \Unip(3,\ol{\varphi})\\
\hline
\HHH_3& \Unip(2)^2&
\varphi \oplus \1\oplus \Unip(1,-{\varphi})^2&
\varphi^2\oplus \Unip(3,{\varphi})\\
\hline 
\HHH_4& \1^2\oplus -\1^3&
-\varphi\oplus \Unip(1,{\varphi})^2\oplus \1^2&
\varphi^2\oplus \1\oplus \Unip(3,{\varphi})\\
\hline 
\HHH_5 &
\Unip(1,\ol{\varphi})^3\oplus \1^3& 
\Unip(2,-\ol{\varphi})\oplus \1^2\oplus \Unip(1,-\ol{\varphi})^2&
\1\oplus \Unip(2,\ol{\varphi}^2)\oplus
\Unip(3,\ol{\varphi})\\
\hline
\HHH_6 & -\1^4\oplus \1^3&
\Unip(2)^2\oplus\Unip(3)&\Unip(2,\varphi)\oplus \Unip(2,\ol{\varphi})
\oplus \Unip(3) \\
\hline
\end{array}$$  
\vspace{.5cm}

\noindent {\it Case 5:} The characters 
$$\varphi, \eta , \varphi\eta , 
\overline{\varphi\eta}, \overline{\eta}, \overline{\varphi}, \1$$
 are pairwise disjoint and 
$$\varphi\ol{\eta}\not= -\1\not= \varphi\eta^2,\varphi^2\eta.$$ Then 
$$ \begin{array}{|c|c|c|c|}
\hline 
&{\rm at}\quad \alpha_1&{\rm at}\quad \alpha_2&{\rm at}\quad \infty\\
\hline
\HHH_0& -\1&-\varphi\eta&{\varphi\eta}\\
\hline 
\HHH_1&  \ol{\varphi\eta^2}\oplus \1 & 
-\overline{\varphi\eta}\oplus -\overline{\varphi}
& \ol{\varphi}\ol{\eta}\oplus \ol{\varphi^2\eta^2}\\
\hline 
\HHH_2&\1\oplus -\1^2 & \eta\oplus \eta^2\oplus \1&\eta\oplus \ol{\varphi}\oplus \varphi\eta^2 \\
\hline
\HHH_3&\Unip(1,\ol{\varphi\eta})^2\oplus \1^2 & 
-\ol{\eta}\oplus \1\oplus \Unip(1,-\varphi\ol{\eta})^2
& \ol{\eta}\oplus \ol{\varphi\eta^2}\oplus \varphi\oplus \ol{\eta}^2\\
\hline 
\HHH_4& \Unip(1,-\1)^3\oplus \1^2 &
-\varphi\oplus \Unip(1,\varphi)^2\oplus \1^2
& \varphi\oplus \ol{\eta}\oplus \varphi^2\eta\oplus \varphi\ol{\eta}
\oplus \varphi\eta\\
\hline 
\HHH_5 &\Unip(1,\ol{\varphi})^3\oplus \1^3 &
\Unip(2,-\ol{\varphi})\oplus \Unip(1,-\ol{\varphi})^2\oplus
 \1^2
 &\ol{\varphi}\oplus \ol{\eta\varphi^2}\oplus \eta \oplus \ol{\varphi\eta}\oplus
\ol{\varphi}\eta\oplus \ol{\varphi}^2 \\
\hline
\HHH_6 & -\1^4\oplus \1^3&
\Unip(2)^2\oplus\Unip(3)&\1\oplus \ol{\eta\varphi} \oplus \eta\varphi\oplus \ol{\eta}\oplus
\eta\oplus \ol{\varphi}\oplus \varphi\\
\hline
\end{array}.$$
By what was said above, 
this finishes the proof of Claim (i).\\

Let us prove Claim (ii):
Let $D$ be a reduced effective divisor on $\AA^1,$ 
let $U:=\AA^1\setminus D,$ and let 
$j:U\to \PP^1$ denote the obvious inclusion. 
A sheaf
$\HHH\in \T_\ell(k)$ which is lisse on $U$ 
is cohomologically rigid, if and only if 
\begin{equation}\label{eqneu0}
{\rm rig}(\HHH)=
(1-\Card(D))\rk(\HHH|_U)^2+\sum_{s\in D\cup \infty}\sum_{i,\chi}e_i(s,\chi,\HHH)^2 =2,
\end{equation} compare to \cite{Katz96}, 6.0.15.
(Note that the sum $\sum_{i,\chi}e_i(s,\chi,\HHH)^2$ gives the 
dimension of the centralizer of the local monodromy in the 
group $\GL_{\rk(\HHH|_U)}(\bQl),$ see loc.~cit. 3.1.15.)
Another necessary condition for $\HHH$ to be contained in $\T_\ell$ is 
that $\HHH$ is irreducible. This implies that 
\begin{equation}\label{tableneu2}
\chi(\PP^1,j_*(\HHH|_{U}))=(1-\Card(D))\rk(\HHH)
+\sum_{s\in D\cup \{\infty\}} e_1(s,\1,\HHH)
\leq 0,\end{equation}
since  the same arguments as on Formula \eqref{eqabove} apply.
Assume first that $\Card(D)>2$ and that 
$\HHH$ fails to be lisse at all points of $D.$
 Then, by \eqref{eqneu0}, 
$$\sum_{s\in D\cup \infty}\sum_{i,\chi}e_i(s,\chi,\HHH)^2=
2+(\Card(D)-1)7^2.$$ Since $\sum_{i,\chi}e_i(s,\chi,\HHH)^2\leq 29$
by Table~1,
one immediately concludes that the cardinality of $D$ is $\leq 3.$ 
If $\Card(D)= 3$ then, by Table~1, the following combinations 
of the centralizer dimensions can occur:
$$(25,25,25,25),\quad (29,29,29,13),\quad (29,29,25,17).$$
In each case one obtains a 
contradiction to \eqref{tableneu2} (using 
a quadratic twist at each local monodromy in the case $(25,25,25,25)$).

Thus 
$D=\{\alpha_1,\alpha_2\},$ where $\alpha_1,\alpha_2$ are two disjoint points  of 
$\AA^1(k),$  and 
$$ \sum_{s\in \{\alpha_1,\alpha_2,\infty\}}\sum_{i,\chi}e_i(s,\chi,\HHH)^2 =7^2+2=51.$$ This leaves one with $7$ possible cases $P_1,\ldots,P_7,$ which
are listed in Table~2.

\begin{table}[ht]
\begin{center}
\begin{tabular}{|c|c|c|c|} 
\hline 
{\rm case}&$\sum_{i,\chi}e_i(\alpha_1,\chi,\HHH)^2$&$\sum_{i,\chi}e_i(\alpha_2,\chi,\HHH)^2$&$\sum_{i,\chi}e_i(\infty,\chi,\HHH)^2$\\
\hline
$P_1$&29&13 &9\\
\hline 
$P_2$&29&11 &11\\
\hline 
$P_3$&25&19 &7\\
\hline 
$P_4$&25&17 &9\\
\hline 
$P_5$&25&13 &13\\
\hline 
$P_6$&19&19 &13\\
\hline 
$P_7$&17&17 &17\\
\hline
\end{tabular}
\vspace{.5cm}

{{\bf Table~2:} {The possible centralizer dimensions}}
%\end{center}
%\caption{The possible centralizer dimensions} \label{T1}
\end{center}
\end{table}
Using Table~1 and the inequality 
\eqref{tableneu2}, one can 
exclude  $P_1,P_4$ and $P_7$ by possibly 
twisting the local monodromy at $\alpha_1,\alpha_2,\infty$ 
by  three suitable (at most quadratic) characters,
whose product is $\1.$ 
The possible case  $P_5$ can be excluded using the inequality in  
 \eqref{tableneu2} and a twist by suitable  characters of order 
at most $4$ whose product is $\1.$  

Since the monodromy representation of $\HHH$ is 
dense in the group $G_2(\bQl),$ one obtains 
an  associated sheaf 
 $\Ad(\HHH)\in \T_\ell$ of generic rank 
$14,$ given by the adjoint representation
of $G_2.$ This is again irreducible, which implies
 that  
\begin{eqnarray*}
\chi(\PP^1,j_*(\Ad(\HHH)|_U))&=&
(1-\Card(D))\cdot 14+
\sum_{s\in D\cup \infty}\dim(C_{G_2}(\HHH(s)))\leq 0.
\end{eqnarray*}  This can be used 
to exclude the case $P_2$ and $P_6,$ because in these cases
one has 
$$ \chi(\PP^1,j_*\underline{\End}_{G_2}(\HHH|_{U}))=-14+8+4+4>0$$
and 
 $$ \chi(\PP^1,j_*\underline{\End}_{G_2}(\HHH|_{U}))=-14+8+8+4> 0 \quad 
\textrm{(resp.).}$$

By the same argument, in case $P_3,$ one can exclude that 
the centralizer dimension $19$ comes from a non-unipotent character. 

Thus, up to a permutation of the points 
$\alpha_1,\alpha_2,\infty,$  we are left with the following possibility 
for the local monodromy: The local monodromy at $\alpha_1$ is an involution,
the local monodromy at $\alpha_2$ 
is unipotent of the form $\Unip(2)^2\oplus \Unip(3)$
and the local monodromy at $\infty$ is regular, i.e., 
the dimension of the centralizer is $7.$ By Table~1,
 Thm.~\ref{thmneuwicht} (i) above 
lists all the possibilities for 
the local monodromy at $\infty$ to be regular, except in in Case 3,
when $\varphi$ has order $6,$ or in Case 4, 
when $\varphi$ is a character of order 4, or in 
Case 5, when $$ \varphi\ol{\eta}= -\1, \quad {\rm or}\quad 
 \varphi\eta^2=-\1,\quad {\rm or} \quad \varphi^2\eta=-\1.$$
  In all these cases one can show that 
no such local system exists by inverting the construction 
of $\HHH$ (using Formula \eqref{eqnneu9} and  middle tensor 
products with the dual sheaves $\F_i^\vee$) 
and deriving a contradiction to \eqref{eqneu0} or \eqref{tableneu2}. 
\Endproof

Let us assume that $\alpha_1=0,\alpha_2=1$ and that 
$k=\bar{\QQ}.$ Let $\iota:\AA^1_{\bar{\QQ}}\to \AA^1_{\QQ}$
be the basechange map. We call a sheaf $\H$ as in 
Thm.~\ref{thmneuwicht} {\it to be defined over $\QQ,$} 
if $\H|_{\AA_{\bar{\QQ}}\setminus \{0,1\}}$ is of the form 
$\iota^*(\frak{H})$ where $\frak{H}$ is 
lisse  on $\AA^1_{{\QQ}}\setminus \{0,1\}.$

\begin{thm}\label{thmneu8} If a  sheaf $\H$ as in Thm.~\ref{thmneuwicht} (i)
is 
defined over $\QQ,$ then the trace 
of the  local monodromy at $\infty$ 
is contained in $\QQ.$ This is the case, if and only if 
the  local monodromy at $\infty$ is of the following form:
$$ \begin{array}{|c|c|}
\hline 
\Unip(7)&\\
\hline 
\Unip(3,\varphi)
\oplus \Unip(3,\overline{\varphi})\oplus \1& \varphi\mbox{ of order 3}\\
\hline
 \Unip(2,\varphi)
\oplus \Unip(2,\ol{\varphi})\oplus \Unip(3)&  \varphi\mbox{ of order 3 or 6}\\
\hline
\varphi\oplus \eta \oplus \varphi \eta \oplus 
\overline{\varphi\eta}\oplus \overline{\eta}\oplus \overline{\varphi}\oplus 
\1 &
\varphi\mbox{ of order 8 and } \eta=\varphi^2\\
&\varphi \mbox{ of order 7 or 14 and } \eta=\varphi^2\\
& \varphi \mbox{ of order 12 and } \eta=-\varphi\\
\hline \end{array}$$
\end{thm}

\proof The first claim follows from the structure of the 
local fundamental group at $\infty.$ Using the local
monodromy of the sheaves $\HHH$ listed in Thm.~\ref{thmneuwicht} (i),
 one obtains the result by an explicit computation. 
\Endproof

\section{The motivic interpretation of the rigid
$G_2$-sheaves}
In this section we recall the motivic interpretation of the 
middle convolution in the universal setup 
of \cite{Katz96}, Chap. 8. This leads to 
an explicit geometric construction of the rigid $G_2$-sheaves 
found in the previous section.

\subsection{Basic definitions}\label{secreview}
Let us recall the setup of \cite{Katz96}, Chap. 8: 
Let $k$ denote an algebraically closed field and $\ell$ a
prime number which is invertible  in $k.$ Let 
further $\alpha_1,\ldots,\alpha_n$ be pairwise disjoint points
of $\AA^1(k)$ and $\zeta$ a primitive root of unity in 
$k.$ 
Fix an integer $N\geq 1$ such that $\chara(k)$ does not divide 
$N$ and 
let
$$ R:=R_{N,\ell}:=\ZZ[\zeta_N,\frac{1}{N\ell}],$$
where $\zeta_N$ denotes a primitive $N$-th root of unity. Set
$$ S_{N,n,\ell}:=R_{N,\ell}[T_1,\ldots,T_n][1/\Delta],\,\,
\Delta:=\prod_{i\not= j}(T_i-T_j).$$
Fix an embedding $R\to \ol{\QQ}_\ell$ and let $E$ denote the fraction
field of $R.$ For a place $\lambda$ of  $E,$ let $E_\lambda$ denote the 
$\lambda$-adic completion of $E.$ Let 
$\phi:S_{N,n,\ell}\to k$ denote the unique ring homomorphism 
for which $\phi(\zeta_N)=\zeta$ and for which
$$\phi(T_i)=\alpha_i,\quad i=1,\ldots,n.$$
Let $\AA^1_{S_{N,n,\ell}}\setminus \{T_1,\ldots,T_n\}$
denote the affine line over $S$ 
with the $n$ sections $T_1,\ldots,T_n$
deleted. 
Consider more generally the spaces 
$$\AA(n,r+1)_{R}:=\Spec(R[T_1,\ldots,T_n,X_1,\ldots,X_{r+1}][\frac{1}{\Delta_{n,r}}]),$$
where
$$ \Delta_{n,r}:=(\prod_{i\not= j}(T_i-T_j))(\prod_{a,j}(X_a-T_j))
(\prod_{k}(X_{k+1}-X_k))$$ (here the indices $i,j$ run 
through $\{1,\ldots,n\},$ the index 
$a$ through $\{1,\ldots,r+1\}$ and the index $k$ runs through $\{1,\ldots,r\};$
when $r=0$ the empty product $\prod_k(X_{k+1}-X_k)$ is understood to be $1$).

Let 
$$\pr_i:\AA(n,r+1)_{R} \to \AA^1_{S_{N,n,\ell}}\setminus \{T_1,\ldots,T_n\},$$
$$ (T_1,\ldots,T_n,X_1,\ldots,X_{r+1})\mapsto (T_1,\ldots,T_n,X_i).$$

On $(\GG_m)_R$ with coordinate $Z,$ one has the Kummer covering 
of degree $N,$ of equation $Y^N=Z.$ This is a connected $\mu_N(R)$-torsor
whose existence defines a surjective 
homomorphism 
$\pi_1((\GG_m)_R)\to \mu_N(R).$ The  chosen embedding 
$R\to \ol{\QQ}_\ell$ defines a faithful
character
$$ \chi_N: \mu_N(R) \to \ol{\QQ}_\ell^\times$$ and the 
composite homomorphism 
$$ \pi_1((\GG_m)_R)\to \mu_N(R) \to \ol{\QQ}_\ell^\times$$
defines the Kummer sheaf $\L_{\chi_N}$ on $(\GG_m)_R.$ For any scheme
$W$ and any map $f: W\to (\GG_m)_R,$ define 
$$\L_{\chi(f)}:=f^*\L_{\chi}.$$

\subsection{The middle convolution of local systems}\label{secreview1}

Denote by $\LisseNnl$ the category of lisse
$\bQl$-sheaves on 
$$ \AA(n,1)_R=(\AA^1-(T_1,\ldots,T_n))_{S_{N,n,\ell}}.$$ For each nontrivial 
$\ol{\QQ}_\ell$-valued character $\chi$ of the group 
$\mu_N(R),$ Katz \cite{Katz96} 
defines a left exact {\it middle convolution functor}
$$ \MC_\chi: \LisseNnl \To \LisseNnl$$
as follows: 

\begin{defn} \label{defdef}{\rm View the space 
$\AA(n,2)_R$ with its second projection 
$\pr_2$ to $\AA(n,1)_R,$ 
as a relative $\AA^1$ with coordinate $X_1,$ minus 
the $n+1$ sections $T_1,\ldots,T_n,X_2.$ Compactify the 
morphism $\pr_2$ into the relative $\PP^1$
$$ \ol{\pr}_2: \PP^1\times \AA(n,1)_R \to \AA(n,1)_R,$$
by filling in the sections $T_1,\ldots,T_n,X_2,\infty.$
Moreover, let $j:\AA(n,2)_R \To \PP^1\times \AA(n,1)_R$ denote the 
natural inclusion. 
The
{\it middle convolution of $\FFF \in \LisseNnl$ and $\L_\chi$} is
defined as follows 
$$ \MC_\chi(\FFF):=R^1(\ol{\pr}_2)_!
(j_*(\pr_1^*(\FFF)\otimes \L_{\chi(X_2-X_1)})) \in \LisseNnl,$$
see loc.~cit. Section~8.3.} \end{defn}

For any $\FFF\in \LisseNnl,$ and any nontrivial 
character $\chi$ as above, let $\FFF_k$ 
 denote the restriction of $\FFF$ to the geometric fibre
$U_k=\AA^1_k\setminus \{\alpha_1,\ldots,\alpha_n\}$ of 
$(\AA^1-(T_1,\ldots,T_n))_{S_{N,n,\ell}}$ which is 
defined by the homomorphism 
$\phi: S\to k.$ Define $\chi_k$ as the restriction of 
$\chi$ to $\GG_{m,k}$ and let  $j:U_k\to \PP^1_k$ denote the inclusion.
Then the following holds:
\begin{equation}\label{eqeqeq} \MC_{\chi_k}(j_*\FFF_k)|_{U_k}=\MC_\chi(\FFF)_k,\end{equation}
where on the left, the middle convolution
$\MC_{\chi_k}(\FFF_k)$ is defined as in Section
\ref{Secneu1} and on the right, the middle convolution is defined
as in Def. \ref{defdef} above (see \cite{Katz96}, Lemma~8.3.2). \\

\subsection{The motivic interpretation of the middle convolution}
In \cite{Katz96}, Thm.~8.3.5 and Thm.~8.4.1,
 the following result is proved:

\begin{thm}\label{thmmotivicint} Fix an integer $r\geq 0.$ For a choice of 
$n(r+1)$ characters 
$$ \chi_{a,i}:\mu_N(R)\to \bQl^\times,\quad i=1,\ldots,n,\quad 
a=1,\ldots,r+1,$$ and a choice 
of $r$ nontrivial characters
$$ \rho_k: \mu_N(R)\to \bQl^\times,\quad k=1,\ldots,r,$$ define a rank one sheaf $\L$ on
$\AA(n,r+1)_R$ by setting
$$ \L:=\bigotimes_{a,i}\L_{\chi_{a,i}(X_a-T_i)}
\bigotimes_{k}\L_{\rho_k(X_{k+1}-X_k)}.$$ Then 
the following holds:
\begin{enumerate}
\item The sheaf $\K:=R^r(\pr_{r+1})_!(\L)$ is mixed of integral weights in 
$[0,r].$ There exists a short exact sequence of lisse sheaves on 
$\AA^1_{S_{N,n,\ell}}\setminus \{T_1,\ldots,T_n\}:$
$$ 0\to \K_{\leq r-1}\to \K\to \K_{=r}\to 0,$$ 
such that $\K_{\leq r-1}$ is mixed of integral weights $\leq r-1$ 
and where $\K_{=r}$ is punctually pure of weight $r.$ 
\item Let $\chi=\chi_N:\mu_N(R)\to \bQl^\times $ be the faithful character
defined in the last section and let 
$e(a,i),  \, i=1,\ldots,n,\, 
a=1,\ldots,r+1,$ and $f(k),\, k=1,\ldots,r$ be integers with 
$$ \chi_{a,i}=\chi^{e(a,i)},\quad \textrm{and}\quad \rho_k=\chi^{f(k)}.$$ 
In the product space $\GG_{m,R}\times \AA(n,r+1)_R,$ consider the
hypersurface ${\rm Hyp}$ given by the equation
$$ Y^N=\left(\prod_{a,i}(X_a-T_i)^{e(a,i)}\right)\left(
\prod_{k=1,\ldots,r} (X_{k+1}-X_k)^{f(k)}\right)$$
and let $$\pi: {\rm Hyp}\to \AA^1_{S_{N,n,\ell}}\setminus \{T_1,\ldots,T_n\},$$
$$ (Y,T_1,\ldots,T_n,X_1,\ldots,X_{r+1})
\mapsto (T_1,\ldots,T_n,X_{r+1}).$$ The group $\mu_N(R)$ 
acts on $\Hyp$ by permuting $Y$ alone, inducing an action 
of  $\mu_N(R)$ on $R^r\pi_!(\bQl).$ Then the sheaf   $\K$ is isomorphic to the 
$\chi$-component 
$ \left(R^r\pi_!\bQl \right)^\chi$ of 
$R^r\pi_!(\bQl).$   
\item  For $a=1,\ldots,r+1,$ 
let 
$$\F_a=\F_{a}(X_a):=\bigotimes_{i=1,\ldots,n}\L_{\chi_{a,i}(X_a-T_i)}\in \LisseNnl.$$ 
Let 
\begin{eqnarray*}\HHH_0&:=&\F_1, \\
\HHH_1&:=&\F_2\otimes \MC_{\rho_1}(\HHH_0),\\
& \vdots&\\
\HHH_r&:=&\F_{r+1}\otimes \MC_{\rho_r}(\HHH_{r-1}).\end{eqnarray*}
 Then 
$ \K_{=r}=\HHH_r.$
\end{enumerate}
\end{thm}

\subsection{Motivic interpretation for rigid $G_2$-sheaves}

Let $\epsilon:\pi_1(\GG_{m,R})\to \mu_N(R)$ be the 
surjective homomorphism of Section~\ref{secreview}. 
By composition with $\epsilon,$ every character
$\chi:\mu_N(R)\to \bQl^\times$ gives rise to a character
of $\pi_1(\GG_{m,R}),$ again denoted by $\chi.$
For a sheaf $\K\in \LisseNnl$ which is mixed of integral weights 
in $[0,r]$
let $W^r(\K)$ denote the weight-$r$-quotient of $\K.$ 

\begin{thm} \label{thmmotg2} Let $\varphi,\eta$ and 
$\HHH(\varphi,\eta)$ be as in Thm.~\ref{thmneuwicht}. Let $N$ denote the 
least common multiple of $2$ and the orders of 
$\varphi,\eta.$ Let further $\chi=\chi_N:\mu_N(R)\to \bQl^\times$
be the character of order $N$ which is defined in Section 
\ref{secreview}, and  
let  $n_1,n_2$ be integers with 
$$\varphi =\chi_k^{n_1}\quad {\rm and} \quad \eta=\chi_k^{n_2},\quad n_1,n_2\in \ZZ,$$ 
where $\chi_k$ is the restriction of $\chi$ to $\GG_{m,k}.$
Let $\Hyp=\Hyp(n_1,n_2)$ denote
the hypersurface in $\GG_{m,R}\times \AA(2,6+1)_R,$ 
given by the following equation:
$$ Y^N=\left(\prod_{ 1\leq a\leq 7;\; 1\leq i\leq 2}(X_a-T_i)^{e(a,i)}\right)\left(
\prod_{1\leq k \leq 6} (X_{k+1}-X_k)^{f(k)}\right),$$
where  the numbers $e(a,i)$ and the $f(k)$ are as follows:
 $$
\begin{array}{|c|c|c|c|c|c|c|}
\hline
\; e(1,1)&e(2,1)&e(3,1)&e(4,1)&e(5,1)&e(6,1)&e(7,1)\;\\
\hline
\quad \frac{N}{2}\quad &0&\quad \frac{N}{2}\quad &0&\quad \frac{N}{2}\quad &0&\quad \frac{N}{2}\quad \\
\hline
\hline
e(1,2)&e(2,2)&e(3,2)&e(4,2)&e(5,2)&e(6,2)&e(7,2)\\
\hline
\frac{N}{2}+n_1+n_2&\frac{N}{2}-n_1&
0&\frac{N}{2}+n_1-n_2&0&\frac{N}{2}-n_1&0\\
\hline
\end{array}$$
and 
$$ \begin{array}{|c|c|c|c|c|c|}
\hline
f(1)&f(2)&f(3)&f(4)&f(5)&f(6)\\
\hline
\frac{N}{2}-n_1-2n_2&\frac{N}{2}+n_1+2n_2&
\frac{N}{2}-n_1-n_2& \frac{N}{2}+n_1+n_2&
\frac{N}{2}-n_1&\frac{N}{2}+n_1\\
\hline
\end{array}\quad .$$
Let  $$\pi=\pi(n_1,n_2): {\rm Hyp}(n_1,n_2)\,\,\to \,\,\AA^1_{S_{N,n,\ell}}\setminus \{T_1,T_2\},$$ given by 
$
(Y,T_1,T_2,X_1,\ldots,X_{7})\mapsto (T_1,T_2,X_{7}).$
 Then the higher direct image sheaf 
$\,W^6\left[\left(R^6\pi_!\bQl\right)^\chi\right]$ is contained 
in $\LisseNnl.$ Moreover, for any algebraically closed field $k$
whose characteristic does not divide $\ell N,$
one has  an isomorphism 
\begin{eqnarray*} \HHH(\varphi,\eta)|_{\AA^1_k\setminus \{\alpha_1,\alpha_2\}}=
\left(W^6\left[\left(R^6\pi_!\bQl\right)^\chi\right]\right)|_{\AA^1_k\setminus \{\alpha_1,\alpha_2\}}.
\end{eqnarray*}
\end{thm}
\proof This is just a restatement of Thm.~\ref{thmmotivicint}
in the situation of Thm.~\ref{thmneuwicht}. The last formula follows
from Formula \eqref{eqeqeq} and Thm.~\ref{thmmotivicint} (iii).
\Endproof

We now turn to the special case, where $n_1=n_2=0,$  $N=2,$
and $\alpha_1=0,\,\alpha_2=1.$  
In this case,  the higher direct image sheaves which occur in 
Thm.~\ref{thmmotg2} can be expressed
in terms of the cohomology of a smooth {\it and proper} map of schemes
over $\QQ.$  
This will be crucial 
in the next section.

\begin{cor} \label{cormotg4} Let 
$N=2,\, n_1=n_2=0,$
and let 
$$\Hyp=\Hyp(0,0)\subseteq  \GG_{m,R}\times \AA(2,6+1)_R$$
be the associated hypersurface equipped with the structural 
morphism $\pi=\pi(0,0):\Hyp\to \AA^1_{S}\setminus \{T_1,T_2\}.$
Let $\pi_\QQ:\Hyp_\QQ\to \AA^1_{\QQ}\setminus \{0,1\},$
denote the basechange of $\pi$ induced by
$ T_1\mapsto 0$ and $T_2\mapsto 1.$
Then the following holds:
\begin{enumerate}
\item There 
exists a smooth and projective scheme $X$  over 
$\AA^1_{\QQ}\setminus \{0,1\}$ and  an open embedding of 
$j:\Hyp_\QQ \to X$
 such that 
$$D=X\setminus \Hyp_\QQ=\bigcup_{i\in I}D_i$$ 
is a strict
normal crossings divisor over $\AA^1\setminus \{0,1\}.$ 
The involutory 
automorphism $\sigma$ of $\Hyp$ (given by $Y\mapsto -Y$)
 extends to an automorphism $\sigma$ of 
$X.$
 \item Let $\coprod_{i\in I}D_i$ denote the disjoint union 
of the components of $D$ and let 
$$\pi_X:X\to \AA^1_{\QQ}\setminus \{0,1\}\quad{\rm and}\quad 
 \pi_{\coprod D_i}: \coprod_{i\in I}D_i\to \AA^1_{\QQ}\setminus \{0,1\}$$ 
denote the structural 
morphisms. Let $\GGG:=W^6\left[\left(R^6\pi_!\bQl\right)^\chi\right]|_{\AA^1_\QQ\setminus \{0,1\}}.$ 
Then 
$$ \GGG \simeq  \Pi\left[\ker\left(R^6(\pi_{X})_*(\bQl) \to R^6(\pi_{\coprod D_i})_*(\bQl)\right)\right],$$
where $\Pi$ denotes the formal sum $\frac{1}{2}(\sigma - 1).$ 
\end{enumerate}
\end{cor}

\proof Let $\Delta\subseteq \AA^7$ be the divisor defined by the
vanishing of 
\begin{equation}\label{laebl0}
\prod_{i=1}^6 (X_{i+1}-X_i)\prod_{i=1}^7X_i\prod_{i=1}^7
(X_i-1).\end{equation}
By Thm.~\ref{thmmotg2}, 
the hypersurface $\Hyp_\QQ$ is an unramified  double cover of 
an open subset of $\AA^7\setminus \Delta$ defined by 
\begin{equation}\label{laebl} Y^2=\prod_{i=1}^6 (X_{i+1}-X_i)\prod_{i=1,3,5,7}X_i\prod_{i=1,2,4,6}
(X_i-1).\end{equation}
This defines a ramified 
double cover  $\alpha:\ol{X}\to  \PP^6_S=\PP^6\times S,$ where 
$S:=\AA^1_{X_7}\setminus \{0,1\}.$ 
The image of the complement $\ol{X}\setminus \Hyp$ under $\alpha$ 
is a relative divisor $L$ over $S$
on $\PP^6_S.$ The divisor $L$  is  the union of the relative 
hyperplane at infinity 
$L_0=\PP^6_S\setminus (\AA^6_{X_1,\ldots,X_6}\times S)$ with the  $20$ linear hyperplanes 
$L_i=L_{T_i},\, i=1,\ldots ,20,$
which are defined by the vanishing of the partial  projectivation of the irreducible 
factors $T_i$ of  the right hand side
of Equation \eqref{laebl0}. The singularities of $\ol{X}$ are 
situated over the singularities of the ramification locus $R$ of 
$\alpha:\ol{X}\to \PP^6_S$
 which is a subdivisor  of $L$ by  Equation \eqref{laebl}.

There is a standard resolution of any linear 
hyperplane arrangement $L=\bigcup_i L_i\subseteq \PP^n$ 
given in 
\cite{ESV}, Section~2. By this we mean a birational map $\tau:\tilde{\PP}^n
\to \PP^n$ 
which factors into several blow ups and  which has the following properties:
The inverse image
of  $L$ under $\tau$ 
is a strict normal crossings divisor  in $\tilde{\PP}^n$ and the strict transform of $L$ is nonsingular 
(see \cite{ESV}, Claim in Section~2). 
 The standard resolution depends only on the combinatorial 
intersection behaviour of the irreducible components $L_i$ of $L,$ therefore
it can be defined for locally trivial families of hyperplane 
arrangements. 

In our case, we obtain
 a birational map 
  $\tau:\tilde{\PP}^6_S\to \PP^6_S$ 
such that $\tilde{L}:=\tau^{-1}(L)$ 
is a relative strict normal crossings divisor
over $S$ and such that the strict transform of $L$ is smooth over $S.$  
Let $\tilde{\alpha}:\tilde{X}\to \tilde{\PP}^6_S$ 
denote the pullback of the double cover $\alpha$ along
$\tau$ and let $\tilde{R}$ be the ramification divisor of $\tilde{\alpha}.$ 
Then
$\tilde{R}$  is 
a relative strict  normal crossings divisor since it is 
contained in $\tilde{L}.$ 
  Write $\tilde{R}$ as a union 
$\bigcup_k\tilde{R}_k$ of irreducible components. By
 successively 
blowing up the (strict transforms of the) intersection loci $\tilde{R}_{k_1}\cap \tilde{R}_{k_2},\, k_1<k_2,$ one 
ends up with a birational map $\hat{f}: \hat{\PP}^6_S\to \tilde{\PP}^6_S.$
Let $\hat{\alpha}:X\to  \hat{\PP}^6_S$ denote the pullback of the double cover 
$\tilde{\alpha}$ along $\hat{f}.$  Then the strict transform 
of $\tilde{R}$ in $\hat{\PP}^6_S$ 
is a disjoint union of smooth components. Moreover,   
since $\tilde{R}$ is a normal crossings divisor, the exceptional divisor 
of the map $\hat{f}$ has no components 
in common with the ramification locus of $\hat{\alpha}.$
It follows that the double cover $\hat{\alpha}:X\to \hat{\PP}^6_S$ 
is smooth over $S$ and that $D=X\setminus \Hyp$ 
is a strict normal crossings divisor 
over $S.$ This desingularization is
obviously equivariant with respect to $\sigma$ which finishes the proof 
of Claim (i).

Let $\pi_X:X\to S$ denote the structural map 
(the composition of 
$\hat{\alpha}:X\to \hat{\PP}^6_S$ with the natural map $\hat{\PP}^6_S\to S$).
There exists an $n\in \NN$ such that 
the morphism $\pi_X$ extends to a morphism 
$X_A\to \AA^1_{A}\setminus \{0,1\} $
of schemes over $A:=\ZZ[\frac{1}{2n}].$
 We  assume that $n$ is big enough that 
$D_A:=X_A\setminus \Hyp_A$ is a normal crossings divisor over $\AA^1_{A}\setminus \{0,1\}.$ In the following, we mostly omit the subscript ${}_A$ but 
we will tacitly work in the category of schemes over $A$ (making
use of the fact  that $A$ is finitely generated over $\ZZ,$
in order to be able to apply Deligne's results on the Weil conjectures). 
 Let $\pi_{D}:D\to\AA^1\setminus \{0,1\} $ (resp. 
$\pi_{{\Hyp}}:{\Hyp}\to \AA^1\setminus \{0,1\}$) be the 
structural morphisms. 
The excision sequence gives an exact sequence of sheaves
\begin{equation}\label{eq123}
R^5 (\pi_{D})_*(\bQl)\to R^6(\pi_{\Hyp})_!(\bQl)\to R^6(\pi_X)_* (\bQl)\to R^6(\pi_D)_*(\bQl)\to R^7(\pi_{\Hyp})_!(\bQl).
\end{equation}
By exactness and the work of Deligne (Weil II, \cite{DeligneWeil2}),
the kernel of the map $R^6(\pi_{\Hyp})_!(\bQl)\to R^6(\pi_X)_* (\bQl)$ is 
an integral constructible sheaf which is  
mixed of weights $\leq 5.$ Thus \eqref{eq123} implies an isomorphism 
\begin{equation}\label{eq124} W^6(R^6(\pi_{\Hyp})_!(\bQl)) \To \im\left(R^6(\pi_{\Hyp})_!(\bQl)\to R^6(\pi_X)_* (\bQl)\right).
\end{equation} By the exactness of \eqref{eq123} and by functoriality, one thus 
obtains the following chain of isomorphisms
\begin{eqnarray}\label{eq125} W^6(R^6(\pi_{\Hyp})_!(\bQl))^\chi&\simeq & \im(R^6(\pi_{\Hyp})_!(\bQl)\to R^6(\pi_X)_* (\bQl))^\chi\\
&\simeq& \ker(R^6(\pi_X)_* (\bQl)\to R^6(\pi_D)_*(\bQl))^\chi,\nonumber 
\end{eqnarray}
where the superscript ${}^\chi$ stands for the $\chi$-component of the higher direct 
image in the sense of  \ref{thmmotivicint} (the notion extends in an obvious way to 
$X$ and to $D$).

We claim that the natural map 
\begin{equation}\label{eq126}
\ker(R^6(\pi_X)_* (\bQl)\to R^6(\pi_D)_*(\bQl))^\chi \To \ker(R^6(\pi_X)_* (\bQl)\to R^6(\pi_{ \coprod_{i}D_{A,i}})_*(\bQl))^\chi.
\end{equation}
is an isomorphism. For this we argue as follows:
Since the sheaf $W^6(R^6(\pi_{\Hyp})_!(\bQl))^\chi$ is lisse
(see Thm.\ref{thmmotg2}), the isomorphisms given in 
 \eqref{eq125} imply that 
$$\ker(R^6(\pi_X)_* (\bQl)\to R^6(\pi_D)_*(\bQl))^\chi$$ is lisse. It follows from 
proper base change that 
$$\ker(R^6(\pi_X)_* (\bQl)\to R^6(\pi_{ \coprod_{i}D_{A,i}})_*(\bQl))^\chi$$ is lisse. 
Thus, by the 
specialization theorem
(see \cite{Katz90}, 8.18.2), in order to prove that the map in \eqref{eq126}
is an isomorphism, it suffices to show this for any 
{\it closed} geometric point $\bs$ of ${\Hyp}.$ 
In view of \eqref{eq125}, we have thus 
 to show 
that 
\begin{eqnarray}\label{eq127} W^6(H^6_c(\Hyp_\bs,\bQl))^\chi \simeq \ker(H^6(X_\bs,\bQl)\to H^6(\coprod_{i}D_{\bs,i},\bQl))^\chi.\end{eqnarray}
Let $X^0_\bs=X_\bs,$ and for positive natural numbers
$i,$ let $X^i_\bs$ denote the disjoint union of the 
irreducible components of the locus, where $i$ pairwise 
different components of $D_\bs$ meet. It follows from the 
Weil conjectures \cite{DeligneWeil1} that 
the spectral sequence $E_1=H^j(X^i_\bs,\bQl)_\bs\Rightarrow H^{i+j}_c
(U_\bs,\bQl)$ degenerates at $E_2.$ Consequently,
\begin{eqnarray*}W^6(H^{6}_c
({\Hyp}_\bs,\bQl))&\simeq& \ker(H^6(X_\bs,\bQl)\to H^6(\coprod_{i}D_{\bs,i},\bQl)).\end{eqnarray*} This implies \eqref{eq127} and thus proves that the 
map in \eqref{eq126} is an isomorphism as claimed.
So,
\begin{eqnarray*} W^6(R^6\pi_{{\Hyp_A}*}(\bQl))^\chi
 &\simeq&
\ker\left(R^6(\pi_{X_A})_*(\bQl) \to R^6\pi_{\coprod_{i\in I}D_{A,i}*}(\bQl)\right)^\chi\\
&=& \Pi\left(\ker(R^6(\pi_{X_A})_*(\bQl) \to R^6(\pi_{\coprod_{i\in I}D_{A,i}})_*(\bQl))\right),\end{eqnarray*} where the last equality
is a tautology using the representation theory 
of finite (cyclic) groups.
It follows  that 
$$\GGG=W^6(R^6(\pi_\Hyp)_!(\bQl))^\chi|_{\AA^1_\QQ\setminus \{0,1\}} \simeq  
\Pi\left(\ker(R^6(\pi_{X})_*(\bQl) \to R^6(\pi_{\coprod D_i})_*(\bQl))\,\right),$$
as claimed.  
\Endproof

\section{Relative motives with motivic Galois group $G_2$}\label{secmotives}

\subsection{Preliminaries on motives}
For an introduction to the theory of motives, as well as basic properties
and definitions, we refer the reader to the book of Y.~Andr\'e
\cite{Andre06}.  Let $K$ and $E$  denote
a fields of characteristic zero.  
Let $\V_K$ denote
the category of smooth and projective varieties over $K.$
If $X\in \V_K$ is purely $d$-dimensional, denote by 
$\Corr^0(X,X)_E$ the $E$-algebra of  codimension-$d$-cycles
in $X\times X$ modulo homological equivalence 
(the multiplication is given by the usual composition 
of correspondences).  This notion extends by additivity to an arbitrary
object $X\in \V_K.$ 
A {\it Grothendieck motive with values in $E$}  is then a triple 
$M=(X,p,m),$ where $X\in \V_K,$ $m\in \ZZ,$ and where 
$p\in \Corr^0(X,X)_E$ is idempotent. 
For any $X\in \V_K$ one has associated 
 a motive
$h(X)=(X,\Delta(X),0)$  (called the {\it motive of} $X$), where $\Delta(X)\subseteq X\times X$ denotes
the diagonal.

One also has the  theory of {\it motives 
for motivated cycles} due to Y. Andr\'e \cite{Andre96}, where 
the ring of correspondences $\Corr^0(X,X)_E$ is replaced by a 
 larger ring $\Corr^0_\mot(X,X)_E$ of {\it motivated cycles}
by adjoining a certain homological cycle (the Lefschetz
involution) to $\Corr^0(X,X)_E$
(see \cite{Andre96} and  \cite{Andre06}). 
The formal definition
of a motivated cycle is as follows: For $X,Y\in \V_K,$ let $\pr^{XY}_X$ denote the
projection $X\times Y\to X.$ Then a {\em motivated cycle} is an elements
$(\pr^{XY}_Y)_*(\alpha \cup *_{XY}(\beta)) \in H^*(X),$ where $\alpha,\beta$
are $E$-linear combinations of 
algebraic cycles on $X\times Y$ and $*_{XY}$ is the Lefschetz involution
on $H^*(X\times Y)$ 
relatively to the line bundle $\eta_{X\times Y}=[X]\otimes \eta_Y+\eta_X
\otimes [Y]$ (with $\eta_X,$ resp. $\eta_X,$ arbitrary 
ample line bundles on $X,$ resp. 
$Y$). Define $\Corr^0_\mot(X,X)_E$ as the ring of the motivated codimension-$d$-cycles in analogy to $\Corr^0(X,X)_E.$ 
A 
motive for motivated cycles with values in $E$ 
is then a triple  $M=(X,p,m),$
 where $X\in \V_K,$ $m\in \ZZ,$ and where 
$p\in \Corr^0_\mot(X,X)_E$ is idempotent with respect to the composition 
of motivated cycles. 

The category of motivated cycles is a neutral Tannakian category
(\cite{Andre96}, Section~4).  Thus, by the Tannakian formalism (see 
\cite{DeligneTannaka}), every motive 
for motivated cycles $M$ with values in $E$  has 
attached an algebraic group $G_M$ over $E$ to it,
called the {\it motivic Galois group} of $M.$
Similarly, granting Grothendieck's standard conjectures, 
the category of motives has the structure of a Tannakian category.
 Thus, by the Tannakian formalism 
and by assuming the standard 
conjectures, every motive in the Grothendieck sense
$M$ has attached an algebraic group $\tilde{G}_M$ to it,
called the {\it motivic Galois group} of $M.$
The following lemma and the remark following it 
were communicated to the authors by Y.~Andr\'e:

\begin{lem} \label{lemnew1} Let $M=(X,p,n)$ be a motive for motivated cycles
with motivic Galois group $G_M.$ Assume that 
that Grothendieck's standard conjectures 
hold. Then  
the motive $M$ is defined by algebraic cycles
and
the motivic Galois group $\tilde{G}_M$ in the Grothendieck sense coincides
with  the motivic Galois group $G_M$ of motives for motivated cycles.
\end{lem}

\proof The first claim follows from the fact that 
 the standard conjectures predict the algebraicity of the 
Lefschetz involution in the auxiliary spaces 
$X\times X\times Y$ which are used to define the projector $p$ (see \cite{Saavedra}). 
The last claim follows from the following interpretation
of $G_M$ (resp. $\tilde{G}_M$): The motivic Galois group for motivated 
cycles $G_M$ is the stabilizer of all motivated cycles which appear
in the realizations of all submotives of the mixed
tensors  $M^{\otimes n}\otimes (M^*)^{\otimes n},$
where  $M^*$ denotes the dual of $M$ (this can be seen 
using the arguments in \cite{Andre06}, Chap. 6.3).
 Similarly, under the assumption of the standard conjectures,
  the motivic Galois group $\tilde{G}_M$ is the stabilizer of all
algebraic cycles which appear in the realizations of submotives 
of the mixed
tensors of $M,$ see 
\cite{Andre06}, Chap. 6.3.
 Under the standard conjectures these spaces coincide,
so  $G_M=\tilde{G}_M.$
\Endproof

\begin{rem} \label{remextend} The above lemma
 can be strengthened or expanded as follows.
 It is possible to define  
unconditionally and purely in terms of algebraic cycles a group  
which, under the standard conjectures, will indeed be the motivic  
Galois group of the motive $X=(X,{\rm Id},0),$ where 
$X$ is a smooth projective variety. Namely, let 
  $G^{\rm alg}_{X}$ be
the closed subgroup of $\prod_i \GL (H^i(X)) \times   
\GG_m$  which fixes the classes of algebraic cycles on powers of  $X$   
(viewed as elements of $H(X)^{\otimes n}\otimes \QQ(r)$, the factor  
$\GG_m$ acting on $\QQ(1)$ by homotheties). Then the motivic Galois group  
$G_X$ is related to  $G^{\rm alg}_X$ as follows (cf.~\cite{Andre06},~9.1.3):
   $ G_X =  \im( G^{\rm alg}_{X\times Y} \to G^{\rm alg}_X)$ for a  
suitable projective smooth variety $Y$.
Under the standard conjectures, one may take  $Y$ to be a point.
\end{rem}

\subsection{Results on families of motives}\label{secresultsfam}
It is often useful to consider variations of motives over a base.
Suppose one has given the following data:
\begin{enumerate}
\item A smooth and geometrically connected variety $S$ 
over a field $K\subseteq \CC.$ 
\item Smooth and projective $S$-schemes $X$ and $Y$ of relative dimensions
$d_X$ and $d_Y,$ equipped with invertible ample line bundles
$L_X$ and $L_Y.$ 
\item Two $\QQ$-linear combinations $Z_1$ and $Z_2$ of  
of integral codimension-$d_X+d_Y$-subvarieties in $X\times_S X\times_S Y$  
 which are flat over $S$ such that 
the following holds for one (and thus for all) $s\in S(\CC):$ 
The class 
$$q_s:=(\pr_{X_s\times X_s}^{X_s\times X_s \times Y_s})_*\left([(Z_1)_s
\cup *[(X_2)_s]\right) \in H^{2d_X}(X_s\times X_s)(d_x)
\subseteq \End(H^{*}(X_s))$$ 
satisfies $q_s\circ q_s=q_s,$ where $*$  denotes the Lefschetz involution 
relative to $[(L_{X/S})_s]\otimes [Y_s]+[X_s]\otimes [(L_{Y/S})_s].$
\item an integer $j.$
\end{enumerate}
Then the assignment $s\mapsto (X_s,p_s,j),\, s\in S(\CC),$ defines a 
{\it family of motives} in the sense of  \cite{Andre96}, Section~5.2.
The following result is due to Y. Andr\'e
(see \cite{Andre96}, Thm.~5.2 and Section~5.3):

\begin{thm}\label{thmandre} Let $s\mapsto (X_s,p_s,j),\, s\in S(\CC),$ 
be a family of motives with coefficients
in $E$ and  let $H_E(M_s):=p_s(H^*_B(X_s,E))$ denote the 
{\rm $E$-realization} of $M_s,$ where  $H^*_B(X_s,E)$ denotes 
the singular cohomology ring of $X_s(\CC).$ 
Then there exists a meager 
subset
$\Exc\subseteq S(\CC)$ and a local system of algebraic groups
$G_s\leq \Aut(H_E(M_s))$ 
on $S(\CC),$ such that the following holds:
\begin{enumerate}
\item $G_{M_s}\subseteq G_s$ for all $s\in S(\CC).$ 
\item $G_{M_s}=G_s,$ if and only if $s \notin \Exc.$ 
\item $G_s$ contains the image of a subgroup of finite index 
of $\pi_1^{\rm top}(S(\CC),s).$ 
\item Let
$\fg_s'$ denote the Lie algebra of $G_{s},$ and 
let $\frak{h}_s$ denote the Lie algebra of the Zariski closure
of the image of $\pi_1^{\rm top}(S(\CC),s).$ Then
the Lie algebra 
 $\frak{h}_s$ is an ideal in $\fg_s'.$ 
\end{enumerate}
Moreover, if $S$ is an open subscheme of $\PP^n$ which is defined 
over a number field $K,$   then 
$\Exc\cap \PP^n(K)$ is a {\em thin} subset of $\PP^n(K)$ 
 (thin in the sense of \cite{SerreMordellWeil}). \end{thm}

\subsection{Motives with motivic Galois group $G_2.$}

Let us call an algebraic group $G$ which is defined over
a subfield of $\bQl$ 
 to be {\it of type $G_2$} if 
the group of $\bQl$-points $G(\bQl)$ is isomorphic to the simple exceptional algebraic 
group $G_2(\bQl)$ (see \cite{Borel} for the definition of 
the algebraic group $G_2$).  
It is the aim of this section to prove the existence
of motives for motivated cycles having a motivic Galois group  of type $G_2.$ 

We start in the situation of Cor. \ref{cormotg4}: Let
$S:= \AA^1_\QQ\setminus \{0,1\},$  let
 $\pi_\QQ:\Hyp_\QQ\to S$
be  as in Cor. \ref{cormotg4}, and let 
$\pi_X:X\to S$ be the strict normal crossings compactification of $\Hyp$ 
given by Cor. \ref{cormotg4}. Let 
\begin{equation}\label{eqx} \GGG\simeq  \Pi\left[\ker\left(R^6(\pi_{X})_*(\bQl) \to R^6(\pi_{\coprod D_i})_*(\bQl)\right)\right]\end{equation} be as in 
Cor.~\ref{cormotg4},
 where 
the $D_i$ are the components of the normal crossings divisor 
 $D=X\setminus\Hyp$ over $S.$ 
We want to use the right hand side of this isomorphism to 
define a family of motives $(N_s)_{s\in S(\CC)}$
for motivated cycles
such that the $\bQl$-realization of $N_s$ coincides naturally with 
the stalk $ \GGG_s$ of $\GGG.$ 
  This is done in three steps:

\begin{itemize}
\item Let $$\psi_{s}^*: H^*(X_s,\bQl) \To H^*(\coprod_i D_{s,i},\bQl)$$ 
be the 
 map
which is  induced by the tautological map 
$ \psi_s:=\coprod_i D_{s,i}\to  X_s.$
Let $\Gamma_{\psi_s}\in \Corr^0_\mot(X_s,\coprod_i D_{s,i})_\QQ$ be the graph of 
$\psi_s.$ Note that $\Gamma_{\psi_s}$ can be seen as 
 a morphism 
of motives $$\Gamma_{\psi_s}=\psi_{s}^*:h(X_s)\To h(\coprod_i D_{s,i}).$$
Since 
the category of motives for motivated cycles
is abelian (see \cite{Andre96}, Section~4), 
there exists a kernel motive 
$$K_s=(X_s,p_s,0),\quad p_s\in \Corr^0_\mot(X_s,X_s)_\QQ,$$ 
of the morphism $\psi_{s}^*$ such that 
$$p_s(H^*(X_s,\bQl))=
\ker\left( H^*(X_s,\bQl) \to H^*(\coprod_i D_{s,i},\bQl)\right).$$

\item The  K\"unneth projector $\pi^6_{X_s}:H^*(X_s)\to H^i(X_s)$
is also contained in  
$\Corr^0_\mot(X_s,X_s)_\QQ$
(see \cite{Andre96}, Prop.~2.2). 

\item Let $\Pi_s$ denote the following projector in 
$\Corr^0_\mot(X_s,X_s)_\QQ:$ 
$$\Pi_s:=\frac{1}{2}(\Delta(X_s)-\Gamma_{\sigma_s}),$$
where $\Delta(X_s)\subseteq X_s\times X_s$ denotes 
the diagonal of $X_s$ and $\Gamma_{\sigma_s}\leq X_s\times X_s$ 
denotes the graph of  $\sigma_s.$ 
By construction, the action of $\Pi_s$ on $H^6(X_s)$ 
is the same as the action (induced by) 
the idempotent $\Pi=\frac{1}{2}(1-\sigma)$ which 
occurs in Cor. \ref{cormotg4}. 
\end{itemize}

By \eqref{eqx} one has 
$$ \GGG_s=  \Pi\left[
\ker\left(H^6(X_s,\bQl)\to H^6(\coprod_i D_{s,i},\bQl)\right)\right]
,\quad \forall  s\in S(\CC).$$ Thus,
by combining the above arguments, one sees that 
the stalk $\GGG_s$ is the $\bQl$-realization $H_\bQl(N_s)$ of the 
motives
$$ N_s:=(X_s,\Pi_s\cdot  p_s \cdot \pi^6_{X_s},0)\quad {\rm with }
\quad \Pi_s\cdot  p_s \cdot \pi^6_{X_s}\in \Corr^0_\mot(X_s,X_s)_\QQ.$$ 
We set 
$$ M_s:=N_s(3)=(X_s,\Pi_s\cdot  p_s \cdot \pi^6_{X_s},3).$$ 

\begin{thm}\label{thm2new} The motives $M_s, \, s\in S(\CC),$ 
form a family of motives such that  for any  $s \in S(\QQ)$
outside a thin set, the  motive 
$M_s$  has
a motivic Galois group  
of type $G_2.$  \end{thm}

\proof That the motives $(N_s)_{s\in S(\CC)}$ 
form a family of 
motives (in the sense of Section~\ref{secresultsfam})
  can be seen by the following arguments:
Let $\Gamma_\sigma\subseteq  X\times_S X$ be the graph of the automorphism 
$\sigma$ and let $\Delta(X) \subseteq X\times_S X$ be the diagonal. 
By Cor.~\ref{cormotg4}, the projectors $\Pi_s$ arise from the 
$\QQ$-linear combination of schemes
 $\frac{1}{2}(\Delta(X)-\Gamma_\sigma)$ 
over  $S$ via base change to 
$s.$ 
 The K\"unneth projector $\pi^6_{X_s}\in 
\Corr^0_\mot(X_s,X_s)$ is invariant under the action of $\pi_1(S).$ 
It follows thus from the theorem of the fixed part as in \cite{Andre96}, Section~5.1, that $\pi^6_{X_s}$  arises from the restriction 
of the K\"unneth projector $\pi^6_{\bar{X}},$ 
where $\bar{X}$ denotes a normal crossings 
compactification over $\QQ$ of the morphism
$\pi_X:X\to S$ (which exists by Hironaka \cite{Hironaka}). 
 By \cite{Andre96}, Prop.~2.2., the projector
$\pi^6_{\bar{X}}$ is a motivated cycle. Since this cycle gives 
rise to the K\"unneth projector $\pi^6_{X_s}$ on one fibre via restriction,
one can use the local triviality of the family $X/S$ to show that  
  the restriction of 
$\pi^6_{\bar{X}}\in \Corr^0(\bar{X}\times\bar{X})$ 
to $X\times_S X$ gives  rise to a family of motives
$(X_s,\pi^6_{X_s},0).$ A similar argument 
applies to the projectors $p_s.$ Therefore the motives
$M_s=N_s(3),\, s\in S(\CC),$ form indeed a family of motives.

Let $\GGG^{\rm an}$ be the local 
system on $S(\CC)$ which is defined by 
 the composition of the  natural map 
$\pi_1^{\rm top}(S(\CC),s)
\to  \pi_1(S,s)$ with the monodromy representation of $\GGG.$ 
By the comparison isomorphism between singular- and \'etale cohomology, 
the local system $\GGG^{\rm an}$ coincides with the local system 
which is defined by the singular $\bQl$-realizations $H_\bQl(N_s)$ 
of the above family 
$(N_s)_{s\in S(\CC)}.$
It follows  from 
Thm.~\ref{thmmotg2} 
that 
$ \GGG|_{\AA^1_{{\CC}}\setminus \{0,1\}}\simeq \HHH(\1,\1)|_{\AA^1_{{\CC}}\setminus \{0,1\}},$
where $\HHH(\1,\1)$ is as in Thm.~\ref{thmneuwicht}. It follows thus from
Thm.~\ref{thmneuwicht} (i)
that the image of $\pi_1^{\rm top}(S(\CC),s)
\leq \pi_1(S,s)$ in 
$ \Aut(H_\bQl(N_s))\simeq\GL_7(\bQl)$ under the monodromy 
map is Zariski dense in the group $G_2(\bQl).$ 
Since the algebraic group $G_2$ is connected, the Zariski 
closure of the image of every 
subgroup of finite index of $\pi_1^{\rm top}(S(\CC))$ 
coincides also with 
$G_2(\bQl).$ 

By Thm.~\ref{thmandre}, (i) and (ii), and since $S$ is open in $\PP^1,$ 
there exists
a local system $(G_s)_{s\in S(\CC)}$ 
of algebraic groups with $G_s\leq \Aut(H_{\bQl}(M_s))$ 
such that the following holds: 
The motivic Galois group $ G_{M_s}$ is contained in $G_s, $ and there exists a thin subset $\Exc\subseteq \QQ$ 
such that if $s\in \QQ\setminus \Exc,$ then 
$G_s=G_{M_s}.$ 
By Thm.~\ref{thmandre} (iii),
$G_s $ contains a subgroup
of finite 
index of the image of $\pi_1(S(\CC),s).$ 
Thus, by what was said above, 
the group $G_s$ contains the group $G_2$ for all $s\in S(\CC).$
Let $\frak{g}_2$ denote the Lie algebra of the group 
$G_2.$ 
Let $\frak{g}'_s$ denote the Lie algebra of the group
$G_s.$ By Thm.~\ref{thmandre} (iv), the Lie algebra $\fg_2$ is 
an ideal of $\frak{g}'_s.$  It follows from this 
and from $N_{\GL_7}(G_2)=\GG_m\times G_2$ (where $\GG_m$ denotes the 
subgroup of scalars of $\GL_7$) that
$G_{N_s}\leq \GG_m\times G_2$ for all $s\in S(\CC).$
The representation $\rho_{N_s}$ 
of $G_{N_s},$ which belongs to the motive $N_s$ under 
the Tannaka correspondence, is therefore a tensor product
$\chi\otimes \rho,$ where $\chi:G_{N_s}\to \GG_m$ is a character and $\rho:G_{N_s}\to 
\GL_7$ has values in $G_2\leq \GL_7.$ Let $A_s$ denote the dual of the 
motive which belongs to $\chi$ under the Tannaka correspondence.
Then $G_{N_s\otimes A_s}=G_2$ for all $s\in \QQ\setminus \Exc.$ 

We claim that for 
$s\in \QQ\setminus \Exc,$ the motive  $A_s$ is the motive $(\Spec(s),{\rm Id},3):$ 
The Galois representation which is associated to the motive 
$N_s$ is equivalent to that of the stalk of  $\GGG$ at 
$s$ (viewed as $\bar{\QQ}$-point) and is therefore pure of weight 
$6.$ By \cite{KisinWortmann}, Thm.~3.1, any 
rank-one motive over $\QQ$ is a Tate twist of an Artin motive. 
Therefore, the $\ell$-adic realization of 
any rank-one motive over $\QQ$ is a power of the cyclotomic character
with a finite character. 
In order that $G_{N_s\otimes A_s}$ is contained in $G_2,$ 
the $\ell$-adic realization of $A_s$ has to be of the form 
$\epsilon\otimes  \chi_\ell^3,$ where $\epsilon$ is of order $\leq 2.$ 
But if the order of $\epsilon$ is equal to $2,$ we derive a contradiction
to Thm.~1 of the Appendix. It follows that  $A_s$ is the motive $(\Spec(s),{\rm Id},3)$ and that the motivic Galois group of 
$M_s=N_s\otimes A_s=N_s(3)$ is of type $G_2.$ 
\Endproof

\begin{rem}\label{remext1} Under the hypothesis of 
the standard conjectures, Thm.~\ref{thm2new} and Lemma~\ref{lemnew1} 
imply the 
existence of Grothendieck motives whose motivic Galois group is 
of type $G_2.$ 
 Moreover, it follows from
 Rem.~\ref{remextend} that, independently from the standard conjectures,
  there is a projective smooth variety $X$ over $\QQ$ such that the group
$G^{\rm alg}_X$ (which is defined in Rem.~\ref{remextend}) 
has a quotient $G_2$.
\end{rem}

\bibliographystyle{plain} \bibliography{$HOME/Biblio/p}

\begin{thebibliography}{10}

\bibitem{BBD}
J.~Bernstein A.A.~Beilinson and P.~Deligne.
\newblock {\em Faisceaux pervers}, volume 100 of {\em Ast\'erisque}.
\newblock Soc. Math. France, Paris, 1982.

\bibitem{Andre96}
Y.\ Andr\'e.
\newblock Pour une th\'eorie inconditionelle des motifs.
\newblock {\em Publ. Math. IHES}, 83:\ 5--49, 1996.

\bibitem{Andre06}
Y.~Andr\'e.
\newblock {\em Une introduction aux motifs}.
\newblock Panoramas et Synth\`eses 17. Soc. Math. de France, 2004.

\bibitem{Asch}
M.~Aschbacher.
\newblock Chevalley groups of type ${G}_2$ as the group of a trilinear form.
\newblock {\em J. Algebra}, 109:193--259, 1987.

\bibitem{Beukers-Heckman}
F.~Beukers and G.~Heckman.
\newblock Monodromy for the hypergeometric function $\sb n{F}\sb{n-1}$.
\newblock {\em Invent. Math.}, 95(2):\ 325--354, 1989.

\bibitem{Borel}
A.~Borel.
\newblock {\em Linear algebraic groups}, volume 126 of {\em Graduate Texts in
  Mathematics}.
\newblock Springer, New York, 1991.

\bibitem{Chang}
B.~Chang.
\newblock The conjugate classes of {C}hevalley groups of type {$G\sb{2}$}.
\newblock {\em J. Algebra}, 9:4190--211, 1968.

\bibitem{Deligne70}
P.\ Deligne.
\newblock {\em Equations Differ{\'e}ntielles {\`a} Points Singuliers
  R{\'e}guliers}.
\newblock Lecture Notes in Mathematics 163. Springer-Verlag, 1970.

\bibitem{DeligneWeil1}
P.~Deligne.
\newblock La conjecture de {W}eil. {I}.
\newblock {\em Publ. Math. IHES}, 43:\ 273--307, 1974.

\bibitem{DeligneShimura}
P.~Deligne.
\newblock Vari\'et\'es de {S}himura: interpr\'etation modulaire, et techniques
  de construction de mod\`eles canoniques.
\newblock In M.~Aschbacher et~al., editor, {\em Automorphic forms,
  representations and $L$-functions}, volume Proc. Sympos. Pure Math., XXXIII,,
  pages \ 247--289. Amer. Math. Soc., 1977.

\bibitem{DeligneWeil2}
P.~Deligne.
\newblock La conjecture de {W}eil. {II}.
\newblock {\em Publ. Math. IHES}, 52:\ 137--252, 1980.

\bibitem{DeligneTannaka}
P.~Deligne.
\newblock Cat\'egories tannakiennes.
\newblock In {\em The Grothendieck Festschrift}, volume~II of {\em Progr. Math.
  87}, pages \ 111--195. Birkh\"auser, 1990.

\bibitem{DMOS}
P.~Deligne, J.~Milne, A.~Ogus, and K.~Shih.
\newblock {\em {H}odge cycles, motives, and {S}himura varieties}, volume 900 of
  {\em Lecture Notes in Mathematics}.
\newblock Springer-Verlag, 1982.

\bibitem{ESV}
H.~Esnault, V.~Schechtman, and E.~Viehweg.
\newblock Cohomology of local systems on the complement of hyperplanes.
\newblock {\em Invent. Math.}, 109:557--561, 1992.

\bibitem{Jannsen}
U.~Jannsen et~al.
\newblock {\em Motives I, II}.
\newblock Number~55 in Proceedings of Symposia in Pure Mathematics. Amer. Math.
  Soc., 1994.

\bibitem{FeitFong}
W.~Feit and P.~Fong.
\newblock Rational rigidity of {$G\sb 2(p)$} for any prime $p>5$.
\newblock In M.~Aschbacher et~al., editor, {\em Proceedings of the Rutgers
  group theory year, 1983--1984}, pages \ 323--326. Cambridge Univ. Press,
  1985.

\bibitem{GrossSavin}
B.H. Gross and G.~Savin.
\newblock Motives with {G}alois group of type {$G\sb 2$}: an exceptional
  theta-correspondence.
\newblock {\em Comp. Math.}, 114:\ 153--217, 1998.

\bibitem{SGA5}
A.\ Grothendieck.
\newblock {\em {SGA V: Cohomologie $\ell$-adique et fonctions $L$}}.
\newblock Number 589 in Lecture Notes in Math. Springer-Verlag, 1977.

\bibitem{Hironaka}
H.~Hironaka.
\newblock {Resolution of singularities of an algebraic variety over a field of
  characteristic zero. I, II}.
\newblock {\em Ann. Math.}, 79:109--203; 205--326, 1964.

\bibitem{Ince56}
E.L.\ Ince.
\newblock {\em Ordinary differential equations}.
\newblock Dover Publications, 1956.

\bibitem{KatzGKM}
N.M.~\ Katz.
\newblock {\em {G}auss Sums, {K}loosterman Sums, and Monodromy Groups}.
\newblock Annals of Mathematics Studies 116. Princeton University Press, 1988.

\bibitem{Katz90}
N.M.\ Katz.
\newblock {\em Exponential sums and differential equations}.
\newblock Annals of Mathematics Studies 124. Princeton University Press, 1990.

\bibitem{Katz96}
N.M.\ Katz.
\newblock {\em Rigid Local Systems}.
\newblock Annals of Mathematics Studies 139. Princeton University Press, 1996.

\bibitem{KisinWortmann}
M.\ Kisin and S.~Wortmann.
\newblock A note on {A}rtin motives.
\newblock {\em Math. Res. Lett.}, 10:\ 375--389, 2003.

\bibitem{Langlands}
R.P. Langlands.
\newblock Automorphic representations, {S}himura varieties, and motives. {E}in
  {M}\"archen.
\newblock In A.~Borel and W.~Casselman, editors, {\em Automorphic forms,
  representations and $L$-functions}, volume~33 of {\em Proc. Symp. Pure
  Math.}, pages \ 205--246. Amer. Math. Soc., 1979.

\bibitem{Lawther}
R.~Lawther.
\newblock Jordan block sizes of unipotent elements in exceptional algebraic
  groups.
\newblock {\em Comm. Algebra}, 23(11):4125--4156, 1995.

\bibitem{Riemann57}
B.~Riemann.
\newblock Beitr\"age zur {T}heorie der durch die {G}auss'sche {R}eihe
  darstellbaren {F}unctionen.
\newblock {\em Abhandlungen der K\"oniglichen Gesellschaft der Wissenschaften
  zu G\"ottingen}, 7, 1857.

\bibitem{Saavedra}
N.~Saavedra Rivano.
\newblock {\em Categories {T}annakiennes}, volume 265 of {\em Lecture Notes in
  Mathematics}.
\newblock Springer-Verlag, Berlin-Heidelberg-New York, 1972.

\bibitem{SerreMordellWeil}
J.-P. Serre.
\newblock {\em Lectures on the Mordell-Weil theorem}, volume~15 of {\em Aspects
  of Mathematics}.
\newblock Vieweg, Wiesbaden, 1989.

\bibitem{Ser2}
J.-P.\ Serre.
\newblock Propri\'et\'es conjecturales des groupes de {G}alois motiviques et
  des repr\'esentations $\ell$-adiques.
\newblock In U.~Jannsen et~al., editor, {\em Motives I}, number~55 in
  Proceedings of Symposia in Pure Mathematics, pages 377--400, 1994.

\bibitem{ThompsonG2}
J.G. Thompson.
\newblock Rational rigidity of {$G\sb 2(5)$}.
\newblock In M.~Aschbacher et~al., editor, {\em Proceedings of the Rutgers
  group theory year, 1983--1984}, pages \ 321--322. Cambridge Univ. Press,
  1985.

\end{thebibliography}

\end{document}